\newcommand{\et}{\hspace{-0.08in}{\bf .}\hspace{0.1in}}
\newcommand{\BOX}{\hbox {$\sqcap$ \kern -1em $\sqcup$}}
\newcommand{\qed}{\hskip 2em \hbox{\BOX} \vskip 2ex}
\newcommand{\im}{{\rm im}}

\newcommand{\To}{\Rightarrow}

\newcommand{\Cat}{{\rm Cat}}
\newcommand{\LieGrp}{{\rm LieGrp}}
\newcommand{\LieAlg}{{\rm LieAlg}}

\newcommand{\Diff}{{\rm Diff}}
\newcommand{\Vect}{{\rm Vect}}
\newcommand{\Term}{{\rm Term}}
\renewcommand{\to}{\rightarrow}
\newcommand{\tensor}{\otimes}

\newcommand{\maps}{\colon}

\newcommand{\ad}{{\rm ad}}
\newcommand{\tr}{{\rm tr}}
\renewcommand{\deg}{{\rm deg}}

\newcommand{\g}{{\mathfrak g}}

\newcommand{\R}{{\mathbb R}}

\documentclass{article}
\usepackage{amsmath}
\usepackage{euscript}
\usepackage{latexsym}
\usepackage[all, knot]{xy}
\usepackage{amsfonts,amssymb}

\renewcommand{\bar}{\vec}

\newtheorem{theorem}{Theorem}
\newtheorem{defn}[theorem]{Definition}
\newtheorem{lemma}[theorem]{Lemma}
\newtheorem{corollary}[theorem]{Corollary}
\newtheorem{proposition}[theorem]{Proposition}
\newtheorem{example}[theorem]{Example}

\title{Higher-Dimensional Algebra VI: Lie $2$-Algebras}

\author{John C.\ Baez \\
Department of Mathematics,  University of California\\
Riverside, California 92521 \\
USA \\
\\
Alissa S.\ Crans \\
Department of Mathematics,
Loyola Marymount University \\
1 LMU Drive, Suite 2700 \\
Los Angeles, CA 90045 \\
USA \\
\\
email: baez@math.ucr.edu, acrans@lmu.edu }

\date{October 4, 2004}
\begin{document}
\bibliographystyle{plain}
\maketitle

\begin{abstract}
\noindent
The theory of Lie algebras can be categorified starting from a new
notion of `2-vector space', which we define as an internal category in
$\Vect$.  There is a 2-category $2\Vect$ having these 2-vector spaces
as objects, `linear functors' as morphisms and `linear natural
transformations' as 2-morphisms.  We define a `semistrict Lie
2-algebra' to be a 2-vector space $L$ equipped with a skew-symmetric
bilinear functor $[\cdot,\cdot] \maps L \times L \to L$ satisfying
the Jacobi identity up to a completely antisymmetric trilinear
natural transformation called the `Jacobiator', which in turn must
satisfy a certain law of its own.  This law is closely related to
the Zamolodchikov tetrahedron equation, and indeed we prove that
any semistrict Lie 2-algebra gives a solution of this equation,
just as any Lie algebra gives a solution of the Yang--Baxter
equation.  We construct a 2-category of semistrict Lie 2-algebras
and prove that it is 2-equivalent to the 2-category of 2-term
$L_\infty$-algebras in the sense of Stasheff.  We also study
strict and skeletal Lie 2-algebras, obtaining the former from
strict Lie 2-groups and using the latter to classify Lie
2-algebras in terms of 3rd cohomology classes in Lie algebra
cohomology.  This classification allows us to construct for any
finite-dimensional Lie algebra $\g$ a canonical 1-parameter
family of Lie 2-algebras $\g_\hbar$ which reduces to $\g$ at
$\hbar = 0$.  These are closely related to the 2-groups $G_\hbar$
constructed in a companion paper.
\end{abstract}

\section{Introduction}

One of the goals of higher-dimensional algebra is to `categorify'
mathematical concepts, replacing equational laws by isomorphisms
satisfying new coherence laws of their own.  By iterating this process,
we hope to find $n$-categorical and eventually $\omega$-categorical
generalizations of as many mathematical concepts as possible, and
use these to strengthen --- and often simplify --- the connections
between different parts of mathematics.  The previous paper of
this series, HDA5 \cite{BLau}, categorified the concept of Lie group
and began to explore the resulting theory of `Lie 2-groups'.  Here
we do the same for the concept of Lie algebra, obtaining a theory
of `Lie 2-algebras'.

In the theory of groups, associativity plays a crucial role.  When we
categorify the theory of groups, this equational law is replaced by an
isomorphism called the associator, which satisfies a new law of its
own called the pentagon equation.  The counterpart of the associative
law in the theory of Lie algebras is the Jacobi identity.  In a `Lie
2-algebra' this is replaced by an isomorphism which we call the
\emph{Jacobiator}.  This isomorphism satisfies an interesting new law
of its own.  As we shall see, this law, like the pentagon equation,
can be traced back to Stasheff's work on homotopy-invariant algebraic
structures --- in this case, his work on $L_\infty$-algebras, also
known as strongly homotopy Lie algebras \cite{LS,SS}.  This
demonstrates yet again the close connection between categorification
and homotopy theory.

To prepare for our work on Lie 2-algebras, we begin in Section
\ref{intcats} by reviewing the theory of internal categories.  This
gives a systematic way to categorify concepts: if $K$ is some category
of algebraic structures, a `category in $K$' will be one of these
structures but with categories taking the role of sets.
Unfortunately, this internalization process only gives a `strict' way
to categorify, in which equations are replaced by identity morphisms.
Nonetheless it can be a useful first step.

In Section \ref{2vs}, we focus on categories in $\Vect$, the category
of vector spaces.  We boldly call these `$2$-vector spaces', despite
the fact that this term is already used to refer to a very different
categorification of the concept of vector space \cite{KV}, for it is
our contention that our 2-vector spaces lead to a more interesting
version of categorified linear algebra than the traditional ones.  For
example, the tangent space at the identity of a Lie 2-group is a
2-vector space of our sort, and this gives a canonical representation
of the Lie 2-group: its `adjoint representation'.  This is contrast to
the phenomenon observed by Barrett and Mackaay \cite{BM}, namely that
Lie 2-groups have few interesting representations on the traditional
sort of 2-vector space.  One reason for the difference is that the
traditional 2-vector spaces do not have a way to `subtract' objects,
while ours do.  This will be especially important for finding examples
of Lie 2-algebras, since we often wish to set $[x,y] = xy - yx$.

At this point we should admit that our 2-vector spaces are far
from novel entities!  In fact, a category in $\Vect$ is secretly
just the same as a 2-term chain complex of vector spaces.  While
the idea behind this correspondence goes back to Grothendieck
\cite{Gr}, and is by now well-known to the cognoscenti, we
describe it carefully in Proposition \ref{1-1vs}, because it is
crucial for relating `categorified linear algebra' to more
familiar ideas from homological algebra.

In Section \ref{definitions} we introduce the key concept
of `semistrict Lie 2-algebra'.  Roughly speaking, this
is a 2-vector space $L$ equipped with a bilinear functor
\[          [\cdot,\cdot] \maps L \times L \to L ,\]
the Lie bracket, that is
skew-symmetric and satisfies the Jacobi identity up to a
completely antisymmetric trilinear
natural isomorphism, the `Jacobiator' --- which in turn is
required to satisfy a law of its own, the `Jacobiator identity'.
Since we do not weaken the equation $[x,y] = -[y,x]$ to an
isomorphism, we do not reach the more general concept of `weak Lie
2-algebra': this remains a task for the future.

At first the Jacobiator identity may seem rather mysterious.  As
one might expect, it relates two ways of using the Jacobiator to
rebracket an expression of the form $[[[w,x],y],z]$, just as the
pentagon equation relates two ways of using the associator to
reparenthesize an expression of the form $(((w \tensor x) \tensor
y) \tensor z)$.  But its detailed form seems complicated and not
particularly memorable.

However, it turns out that the Jacobiator identity is closely
related to the Zamolodchikov tetrahedron equation, familiar from
the theory of 2-knots and braided monoidal 2-categories
\cite{BLan,BN,CS,Crans,KV}.  In Section \ref{topology} we
prove that just as any Lie algebra gives a solution of the Yang-Baxter
equation, every semistrict Lie 2-algebra gives a solution of the
Zamolodchikov tetrahedron equation!  This pattern suggests that
the theory of `Lie $n$-algebras' --- that is, structures like Lie
algebras with $(n-1)$-categories taking the role of sets --- is
deeply related to the theory of $(n-1)$-dimensional manifolds
embedded in $(n+1)$-dimensional space.

In Section \ref{Linftyalgs}, we recall the definition of an
$L_{\infty}$-algebra.  Briefly, this is a chain complex $V$ of
vector spaces equipped with a bilinear skew-symmetric operation
$[\cdot,\cdot] \maps V \times V \to V$ which satisfies the Jacobi
identity up to an infinite tower of chain homotopies.  We
construct a 2-category of `2-term' $L_\infty$-algebras, that is,
those with $V_i = \{0\}$ except for $i = 0,1$. Finally, we show
this 2-category is equivalent to the previously defined 2-category
of semistrict Lie $2$-algebras.

In the next two sections we study \emph{strict} and \emph{skeletal}
Lie $2$-algebras, the former being those where the Jacobi identity
holds `on the nose', while in the latter, isomorphisms exist only
between identical objects. Section \ref{strictlie2algs} consists of an
introduction to strict Lie $2$-algebras and strict Lie $2$-groups,
together with the process for obtaining the strict Lie $2$-algebra of
a strict Lie $2$-group.  Section \ref{skeletallie2algs} begins with
an exposition of Lie algebra cohomology and its relationship to skeletal
Lie $2$-algebras.  We then show that Lie $2$-algebras can
be classified (up to equivalence) in terms of a Lie algebra
$\mathfrak{g}$, a representation of $\mathfrak{g}$ on a vector
space $V$, and an element of the Lie algebra cohomology group
$H^3(\mathfrak{g},V)$.  With the help of this result, we construct
from any finite-dimensional Lie algebra $\g$ a canonical 1-parameter
family of Lie 2-algebras $\g_\hbar$ which reduces to $\g$ at $\hbar = 0$.
This is a new way of deforming a Lie algebra, in which the
Jacobi identity is weakened in a manner that depends on the parameter
$\hbar$.  It is natural to suspect that this deformation is
related to the theory of quantum groups and affine Lie
algebras.  In HDA5, we give evidence for this by using Chern--Simons
theory to construct 2-groups $G_\hbar$ corresponding
to the Lie 2-algebras $\g_\hbar$ when $\hbar$ is an integer.
However, it would be nice to find a more direct link between
quantum groups, affine Lie algebras and the Lie 2-algebras
$\g_\hbar$.

In Section \ref{conclusions}, we conclude with some guesses
about how the work in this paper should fit into a more general
theory of `$n$-groups' and `Lie $n$-algebras'.

{\bf Note:} In all that follows, we denote the composite of morphisms
$f\maps x \rightarrow y$ and $g\maps y \rightarrow z$ as
$fg\maps x \rightarrow z.$   All 2-categories and 2-functors
referred to in this paper are {\it strict}, though sometimes
we include the word `strict' to emphasize this fact.  We denote
vertical composition of 2-morphisms by juxtaposition;
we denote horizontal composition and whiskering by the symbol $\circ$.


\section{Internal Categories} \label{intcats}

In order to create a hybrid of the notions of a vector space
and a category in the next section, we need the concept of
an `internal category' within some category.
The idea is that given a category $K$, we obtain the definition
of a `category in $K$' by expressing the definition
of a usual (small) category completely in terms of commutative diagrams
and then interpreting those diagrams within $K$.  The same
idea allows us to define functors and natural transformations
in $K$, and ultimately to recapitulate most of category theory,
at least if $K$ has properties sufficiently resembling those
of the category of sets.

Internal categories
were introduced by Ehresmann \cite{E} in the 1960s, and by
now they are a standard part of category theory \cite{Bo}.
However, since not all readers may be familiar with them,
for the sake of a self-contained treatment we start with
the basic definitions.

\begin{defn} \et \label{co}  Let $K$ be a category.
An {\bf internal category} or {\bf category in $K$}, say $X$, consists
of:
\begin{itemize}
\item
an {\bf object of objects} $X_{0} \in K,$
\item
an {\bf object of morphisms} $X_{1} \in K,$
\end{itemize}
together with
\begin{itemize}
\item
{\bf source} and {\bf target} morphisms $s,t \maps X_{1}
\rightarrow X_{0},$
\item
a {\bf identity-assigning} morphism $i \maps X_{0} \rightarrow X_{1},$
\item
a {\bf composition} morphism $\circ \maps X_{1} \times _{X_{0}}
X_{1} \rightarrow X_{1}$
\end{itemize}
such that the following diagrams commute, expressing the usual
category laws:
\begin{itemize}
\item laws specifying the source and target of identity morphisms:
\[
\xymatrix{
    X_{0}
      \ar[r]^{i}
      \ar[dr]_{1}
      & X_{1}
      \ar[d]^{s} \\
     & X_{0} }
\hspace{.2in}
\xymatrix{
      X_{0}
      \ar[r]^{i}
      \ar[dr]_{1}
      & X_{1}
      \ar[d]^{t} \\
     & X_{0}}
\]
\item
laws specifying the source and target of composite morphisms:
\[
  \xymatrix{
   X_{1} \times _{X_{0}} X_{1}
     \ar[rr]^{\circ}
     \ar[dd]_{p_{1}}
     && X_{1}
     \ar[dd]^{s} \\ \\
   X_{1}
     \ar[rr]^{s}
     && X_{0} }
     \hspace{.2in}
\xymatrix{
  X_{1} \times_{X_{0}} X_{1}
     \ar[rr]^{\circ}
     \ar[dd]_{p_{2}}
      && X_{1}
     \ar[dd]^{t} \\ \\
      X_{1}
     \ar[rr]^{t}
      && X_{0} }
\]
\item the associative law for composition of morphisms:
\[
\xymatrix{
   X_{1} \times _{X_{0}} X_{1} \times _{X_{0}} X_{1}
     \ar[rr]^{\circ \times_{X_{0}} 1}
     \ar[dd]_{1 \times_{X_{0}} \circ}
      && X_{1} \times_{X_{0}} X_{1}
     \ar[dd]^{\circ} \\ \\
      X_{1} \times _{X_{0}} X_{1}
     \ar[rr]^{\circ}
      && X_{1} }
\]
\item the left and right unit laws for composition of morphisms:
\[
\xymatrix{
   X_{0} \times _{X_{0}} X_{1}
     \ar[r]^{i \times 1}
     \ar[ddr]_{p_2}
      & X_{1} \times _{X_{0}} X_{1}
     \ar[dd]^{\circ}
      & X_{1} \times_{X_{0}} X_{0}
     \ar[l]_{1 \times i}
     \ar[ddl]^{p_1} \\ \\
      & X_{1} }
\]
\end{itemize}
\end{defn}

The pullbacks referred to in the above definition should be clear from the
usual definition of category; for instance, composition is defined on
pairs of morphisms where the target of the first is the source of the
second, so the pullback $X_1 \times_{X_0} X_1$ is defined via the square
\[
  \xymatrix{
   X_{1} \times_{X_{0}} X_{1}
     \ar[rr]^{p_2}
     \ar[dd]_{p_1}
     && X_{1}
     \ar[dd]^{s} \\ \\
   X_{1}
     \ar[rr]^{t}
     && X_{0} }
\]
Notice that inherent to this definition is the
assumption that the pullbacks involved actually exist.  This
holds automatically when the `ambient category' $K$ has finite
limits, but there are some important examples such as $K = \Diff$
where this is not the case.  Throughout this paper, all of the
categories considered have finite limits:

\begin{itemize}
\item Set, \emph{the category whose objects are sets and whose
morphisms are functions.}

\item Vect, \emph{the category whose objects are vector spaces
over the field $k$ and whose morphisms are linear functions.}

\item Grp, \emph{the category
whose objects are groups and whose morphisms are homomorphisms.}

\item Cat, \emph{the category whose objects are small categories and whose
morphisms are functors.}

\item LieGrp, \emph{the category whose objects are Lie groups and whose
morphisms are Lie group homomorphisms.}

\item LieAlg, \emph{the category whose objects are Lie algebras over
the field $k$ and whose morphisms are Lie algebra homomorphisms.}
\end{itemize}

Having defined `categories in $K$', we can now
internalize the notions of functor and natural transformation
in a similar manner.   We shall use these to construct
a $2$-category $K\Cat$ consisting of categories,
functors, and natural transformations in $K$.

\begin{defn} \et \label{cofunctor} Let $K$ be a category.
Given categories $X$ and $X'$ in $K$, an {\bf internal functor} or
{\bf functor in $K$} between them, say $F\maps X \rightarrow X',$
consists of:
\begin{itemize}
\item
a morphism $F_{0} \maps X_{0} \to X_{0}'$,
\item
a morphism $F_{1} \maps X_{1} \rightarrow X_{1}'$
\end{itemize}
such that the following diagrams commute, corresponding to the
usual laws satisfied by a functor:
\begin{itemize}
\item preservation of source and target:
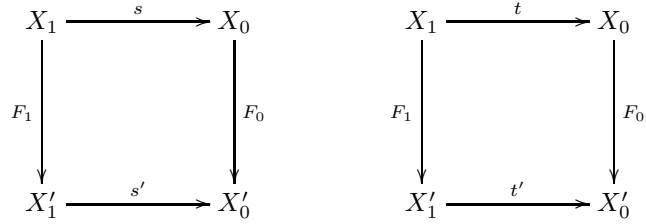
\[
\xymatrix{
  X_{1}
   \ar[rr]^{s}
   \ar[dd]_{F_{1}}
    && X_{0}
   \ar[dd]^{F_{0}} \\ \\
    X_{1}'
   \ar[rr]^{s'}
    && X_{0}' }
\qquad \qquad \xymatrix{
  X_{1}
   \ar[rr]^{t}
   \ar[dd]_{F_{1}}
    && X_{0}
   \ar[dd]^{F_{0}} \\ \\
    X_{1}'
   \ar[rr]^{t'}
    && X_{0}' }
\]
\item preservation of identity morphisms:
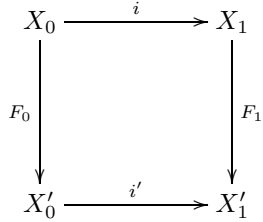
\[
\xymatrix{
    X_{0}
   \ar[rr]^{i}
   \ar[dd]_{F_{0}}
    && X_{1}
   \ar[dd]^{F_{1}} \\ \\
    X_{0}'
   \ar[rr]^{i'}
    && X_{1}' }
\]
\item preservation of composite morphisms:
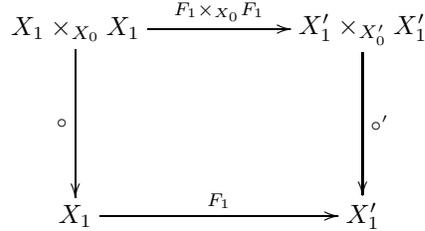
\[
\xymatrix{
   X_{1} \times _{X_{0}} X_{1}
    \ar[rr]^{F_{1} \times_{X_0} F_{1}}
    \ar[dd]_{\circ}
     && X_{1}' \times_{X_{0}'} X_{1}'
    \ar[dd]^{\circ'} \\ \\
     X_{1}
    \ar[rr]^{F_{1}}
     && X_{1}' }
\]
\end{itemize}
\end{defn}

Given two functors $F\maps X \rightarrow X'$ and $G\maps X'
\rightarrow X''$ in some category $K$, we define their composite
$FG\maps X \rightarrow X''$ by taking $(FG)_{0} = F_{0}G_{0}$ and
$(FG)_{1} = F_{1}G_{1}$.  Similarly, we define the identity functor
in $K$, $1_X \maps X \rightarrow X$, by taking $(1_X)_0 = 1_{X_0}$
and $(1_X)_1 = 1_{X_1}$.

\begin{defn} \et \label{conattrans} Let $K$ be a
category.  Given two functors $F,G\maps X \rightarrow X'$ in $K$,
an {\bf internal natural transformation} or {\bf natural
transformation in $K$} between them, say $\theta \maps F
\Rightarrow G$, is a morphism $\theta \maps X_0 \to X'_1$ for
which the following diagrams commute, expressing the usual laws
satisfied by a natural transformation:
\begin{itemize}
\item laws specifying the source and target of a natural transformation:
\[
\xymatrix{
    X_{0}
      \ar[r]^{\theta}
      \ar[dr]_{F_0}
      & X'_{1}
      \ar[d]^{s} \\
     & X_{0} }
\hspace{.2in}
\xymatrix{
      X_{0}
      \ar[r]^{\theta}
      \ar[dr]_{G_0}
      & X'_{1}
      \ar[d]^{t} \\
     & X_{0}}
\]
\item the commutative square law:
\[  \xymatrix{
   X_1
    \ar[rr]^{\Delta (s\theta \times G)}
    \ar[dd]_{\Delta (F \times t\theta)}
     && X'_1 \times_{X_0'} X'_1
    \ar[dd]^{\circ'} \\ \\
     X'_1 \times_{X_0'} X'_1
    \ar[rr]^{\circ'}
     && X'_1
}
\]
\end{itemize}
\end{defn}

Just like ordinary natural transformations, natural
transformations in $K$ may be composed in two different, but
commuting, ways.  First, let $X$ and $X'$ be categories in $K$ and
let $F,G,H\maps X \rightarrow X'$ be functors in $K$.  If $\theta
\maps F \Rightarrow G$ and $\tau \maps G \Rightarrow H$ are
natural transformations in $K$, we define their {\bf vertical}
composite, $\theta \tau \maps F \Rightarrow H,$ by
$$ \theta \tau := \Delta(\theta \times \tau) \circ' .$$
The reader can check that when $K = \Cat$ this reduces to the
usual definition of vertical composition.
We can represent this composite pictorially as:
\[
\xymatrix{
    X
      \ar@/^2pc/[rr]_{\quad}^{F}="1"
      \ar@/_2pc/[rr]_{H}="2"
  && X'
    \ar@{}"1";"2"|(.2){\,}="7"
     \ar@{}"1";"2"|(.8){\,}="8"
    \ar@{=>}"7" ;"8"^{\theta \tau}
  }
\quad = \quad
\xymatrix{
  X
    \ar@/^2pc/[rr]^{F}_{}="0"
    \ar[rr]^<<<<<<{G}^{}="1"
    \ar@/_2pc/[rr]_{H}_{}="2"
    \ar@{=>}"0";"1" ^{\theta}
    \ar@{=>}"1";"2" ^{\tau}
&&  X' }
\]

Next, let $X, X', X''$ be categories in $K$ and let $F,G\maps X
\rightarrow X'$ and $F', G'\maps X' \rightarrow X''$ be functors in
$K$.  If $\theta\maps F \Rightarrow G$ and $\theta'\maps F'
\Rightarrow G'$ are natural transformations in $K$, we define their
{\bf horizontal composite}, $\theta \circ \theta'\maps FF' \Rightarrow
GG'$, in either of two equivalent ways:
\begin{eqnarray*}
  \theta \circ \theta' &:=&
\Delta(F_{0} \times \theta)(\theta' \times G_{1}') \circ'  \\
&=& \Delta(\theta \times G_{0})(F_1' \times \theta') \circ' .
\end{eqnarray*}
Again, this reduces to the usual definition when $K = \Cat$.
The horizontal composite can be depicted as:
\[
\xymatrix{
    X
      \ar@/^2pc/[rr]_{\quad}^{FF'}="1"
      \ar@/_2pc/[rr]_{GG'}="2"
  && X''
    \ar@{}"1";"2"|(.2){\,}="7"
     \ar@{}"1";"2"|(.8){\,}="8"
    \ar@{=>}"7" ;"8"^{\theta \circ \theta '}
  }
\quad = \quad
\xymatrix{
    X
      \ar@/^2pc/[rr]_{\quad}^{F}="1"
      \ar@/_2pc/[rr]_{G}="2"
  && X'
    \ar@{}"1";"2"|(.2){\,}="7"
     \ar@{}"1";"2"|(.8){\,}="8"
    \ar@{=>}"7" ;"8"^{\theta}
    \ar@/^2pc/[rr]_{\quad}^{F'}="1"
    \ar@/_2pc/[rr]_{G'}="2"
  && X''
    \ar@{}"1";"2"|(.2){\,}="7"
     \ar@{}"1";"2"|(.8){\,}="8"
    \ar@{=>}"7" ;"8"^{\theta '}
  }
\]

It is routine to check that these composites are again
natural transformations in $K$.  Finally, given a
functor $F \maps X \rightarrow X'$ in $K$,
the identity natural transformation $1_F \maps F
\Rightarrow F$ in $K$ is given by $1_F = F_0 i$.

We now have all the ingredients of a 2-category:

\begin{proposition}\et Let $K$ be a category.  Then there
exists a strict 2-category \textbf{\textit{K}{\bf Cat}} with
categories in $K$ as objects, functors in $K$ as morphisms,
and natural transformations in $K$ as 2-morphisms, with
composition and identities defined as above.
\end{proposition}

\noindent{\bf Proof. } It is straightforward to check that
all the axioms of a 2-category hold; this result
goes back to Ehresmann \cite{E}. \qed

We now consider internal categories in Vect.


\section{$2$-Vector spaces} \label{2vs}

Since our goal is to categorify the concept of a Lie algebra,
we must first
categorify the concept of a vector space.  A categorified vector
space, or `2-vector space', should be a category with structure
analogous to that of a vector space, with functors replacing
the usual vector space operations.
Kapranov and Voevodsky \cite{KV} implemented this idea
by taking a finite-dimensional 2-vector space to be a category
of the form $\Vect^n$, in analogy to how every finite-dimensional
vector space is of the form $k^n$.
While this idea is useful in contexts such as topological
field theory \cite{Lawrence} and group representation theory \cite{B2},
it has its limitations.  As explained in the Introduction, these
arise from the fact that these 2-vector spaces have
no functor playing the role of `subtraction'.

Here we instead define a 2-vector space to be a category in $\Vect$.
Just as the main ingredient of a Lie algebra is a
vector space, a Lie $2$-algebra will have an underlying $2$-vector
space of this sort.  Thus, in this section we first define a
$2$-category of these $2$-vector spaces.  We then establish the
relationship between these $2$-vector spaces and $2$-term chain
complexes of vector spaces: that is, chain complexes having only
two nonzero vector spaces.  We conclude this section by developing
some `categorified linear algebra' --- the bare minimum necessary
for defining and working with Lie $2$-algebras in the next
section.

In the following we consider vector spaces over an arbitrary
field, $k$.

\begin{defn} \et A {\bf 2-vector space} is a category in {\rm Vect}.
\end{defn}

Thus, a $2$-vector space $V$ is a category with a vector space of
objects $V_0$ and a vector space of morphisms $V_1$, such that the
source and target maps $s,t \maps V_{1} \rightarrow V_{0}$, the
identity-assigning map $i \maps V_{0} \rightarrow V_{1}$, and the
composition map $\circ \maps V_{1} \times _{V_{0}} V_{1}
\rightarrow V_{1}$ are all {\it linear}.  As usual, we write a
morphism as $f \maps x \to y$ when $s(f) = x$ and $t(f) = y$, and
sometimes we write $i(x)$ as $1_x$.

In fact, the structure of a 2-vector space is completely
determined by the vector spaces $V_{0}$ and $V_{1}$
together with the source, target and identity-assigning maps.
As the following lemma demonstrates, composition can always
be expressed in terms of these, together with vector
space addition:

\begin{lemma} \et \label{watereddown} When $K = \Vect$, one can
omit all mention of composition in the definition of category in
$K$, without any effect on the concept being defined.
\end{lemma}

\noindent{\bf Proof. }  First, given vector spaces $V_{0}$, $V_{1}$
and maps $s,t \maps V_{1} \rightarrow V_{0}$ and
$i\maps V_{0} \rightarrow V_{1}$, we will
define a composition operation that satisfies the laws in
Definition \ref{co}, obtaining a 2-vector space.

Given $f \in V_{1}$,
we define the {\bf arrow part} of $f,$ denoted as $\bar{f}$, by
$$\bar{f} = f - i(s(f)).$$
Notice that $\bar{f}$ is in the kernel of the source map since
$$s(f - i(sf)) = s(f) - s(f) = 0.$$
While the source of $\bar{f}$ is always zero, its target may be computed
as follows:
$$t(\bar{f}) = t(f - i(s(f)) = t(f) - s(f).$$
The meaning of the arrow part becomes clearer if we write
$f \maps x \to y$ when $s(f) = x$ and $t(f) = y$.
Then, given any morphism $f \maps x \to y$, we
have $\bar{f} \maps 0 \rightarrow y-x.$
In short, taking the arrow part of $f$ has the effect of `translating
$f$ to the origin'.

We can always recover any morphism from its arrow part together
with its source, since $f = \bar{f} + i(s(f))$.  We shall
take advantage of this by identifying $f \maps x \to y$ with
the ordered pair $(x, \bar{f})$.  Note that with this notation
we have
$$s(x,\bar{f}) = x , \qquad
t(x,\bar{f}) = x + t(\bar{f}).$$

Using this notation,
given morphisms $f \maps x \to y$ and $g\maps y \to z$,
we define their composite by
$$f g := (x, \bar{f} + \bar{g}),$$
or equivalently,
$$(x,\bar{f}) (y,\bar{g}) := (x,\bar{f}+\bar{g}) .$$
It remains to show that with this composition, the
diagrams of Definition \ref{co} commute.  The triangles
specifying the source and target of the identity-assigning morphism
do not involve composition.  The second pair of
diagrams commute since
$$s(f g) = x$$
and
$$t(f g) = x + t(\bar{f}) + t(\bar{g}) = x + (y-x) + (z-y) = z.$$
The associative law holds for composition
because vector space addition is associative.  Finally, the
left unit law is satisfied since given $f \maps x \to y$,
$$i(x) f = (x,0) (x,\bar{f}) = (x,\bar{f}) = f$$
and similarly for the right unit law.  We thus have a 2-vector
space.

Conversely, given a category $V$ in $\Vect$, we shall show that its
composition must be defined by the formula given above.
Suppose that $(f, g)= ((x, \bar{f}), (y, \bar{g}))$ and
$(f', g')= ((x', \bar{f'}), (y', \bar{g'}))$
are composable pairs of morphisms in $V_{1}$.
Since the source and target maps are linear, $(f + f', g + g')$ also
forms a composable pair, and the linearity of composition gives
$$(f + f') (g + g') = f g + f' g'.$$
If we set $g = 1_{y}$ and $f' = 1_{y'}$, the above equation becomes
$$(f + 1_{y'}) (1_y + g') = f 1_y + 1_{y'} g' = f + g'.$$
Expanding out the left hand side we obtain
$$((x, \bar{f}) + (y', 0))
((y,0) + (y', \bar{g'}))
= (x + y', \bar{f}) (y + y', \bar{g'}),$$
while the right hand side becomes
$$(x, \bar{f}) + (y, \bar{g'}) = (x+y', \bar{f} + \bar{g'}).$$
Thus we have $(x + y', \bar{f}) (y + y', \bar{g'})
= (x+y', \bar{f} + \bar{g'})$,
so the formula for composition in an arbitrary $2$-vector space
must be given by
$$f g = (x, \bar{f}) (y, \bar{g}) = (x, \bar{f} + \bar{g})$$
whenever $(f,g)$ is a composable pair.
This shows that we can leave out all reference to composition
in the definition of `category in $K$' without any effect
when $K = \Vect$.
\qed

In order to simplify future arguments, we will often use only the
elements of the above lemma to describe a $2$-vector space.

We continue by defining the morphisms between $2$-vector
spaces:

\begin{defn} \et Given $2$-vector spaces $V$ and
$W$, a {\bf linear functor} $F \maps V \to W$ is a
functor in $\Vect$ from $V$ to $W$.
\end{defn}

\noindent
For now we let $2\Vect$ stand for the category of 2-vector
spaces and linear functors between them; later we will
make $2\Vect$ into a 2-category.

The reader may already have noticed that a 2-vector space
resembles a {\bf 2-term chain complex} of vector spaces: that is,
a pair of vector spaces with a linear map between them,
called the `differential':
$$\xymatrix{
   C_{1} \ar[rr]^{d}
    && C_{0}.  }
$$
In fact, this analogy is very precise. Moreover, it continues
at the level of morphisms. A {\bf chain map} between $2$-term
chain complexes, say $\phi \maps C \to C'$, is simply a pair of
linear maps $\phi_0 \maps C_0 \to C'_0$ and $\phi_1 \maps C_1 \to
C'_1$ that `preserves the differential', meaning that the
following square commutes:
$$\xymatrix{
    C_{1} \ar[rr]^{d}
    \ar[dd]_{\phi_{1}} &&
    C_{0} \ar[dd]^{\phi_{0}} \\ \\
    C_{1}' \ar[rr]^{d'} &&
    C_{0}'
}$$
There is a category $2\Term$ whose objects are 2-term chain complexes
and whose morphisms are chain maps.  Moreover:

\begin{proposition} \et \label{1-1vs}  The categories $2\Vect$
and $2\Term$ are equivalent.
\end{proposition}

\noindent{\bf Proof. } We begin by introducing functors
$$S\maps {\rm 2Vect} \rightarrow {\rm 2Term}$$
and
$$T\maps {\rm 2Term} \rightarrow {\rm 2Vect}.$$
We first define $S$.  Given a 2-vector space $V$, we
define $S(V)=C$ where $C$ is the $2$-term chain complex with
\begin{eqnarray*}
  C_{0} &=& V_{0}, \\
  C_{1} &=& \ker(s) \subseteq V_{1}, \\
  d &=& t|_{C_{1}} ,
\end{eqnarray*}
and $s,t\maps V_{1} \rightarrow V_{0}$ are the source and target
maps associated with the $2$-vector space $V$.  It remains to
define $S$ on morphisms. Let $F \maps V \to V'$ be a linear
functor and let $S(V) = C$, $S(V') = C'.$ We define $S(F)= \phi$
where $\phi$ is the chain map with $\phi_{0} = F_{0}$ and
$\phi_{1} = F_{1}|_{C_{1}}$. Note that $\phi$ preserves the
differential because $F$ preserves the target map.

We now turn to the second functor, $T$.  Given a 2-term
chain complex $C$, we define $T(C) = V$ where $V$ is a
2-vector space with
\begin{eqnarray*}
  V_{0} & = & C_{0}, \\
  V_{1} & = & C_{0} \oplus C_{1}.
\end{eqnarray*}
To completely specify $V$ it suffices by Lemma \ref{watereddown}
to specify linear maps $s,t \maps V_1 \to V_0$ and $i \maps V_0 \to V_1$
and check that $s(i(x)) = t(i(x)) = x$ for all $x \in V_0$.
To define $s$ and $t$, we write any element $f \in V_1$ as a
pair $(x,\bar{f}) \in C_0 \oplus C_1$ and set
\[
\begin{array}{ccccl}
  s(f) &=& s(x, \bar{f}) &=& x, \\
  t(f) &=& t(x, \bar{f}) &=& x + d\bar{f}.
\end{array}
\]
For $i$, we use the same notation and set
\[  i(x) = (x,0)  \]
for all $x \in V_0$.  Clearly $s(i(x)) = t(i(x)) = x$. Note also
that with these definitions, the decomposition $V_1 = C_0 \oplus
C_1$ is precisely the decomposition of morphisms into their source
and `arrow part', as in the proof of Lemma \ref{watereddown}.
Moreover, given any morphism $f = (x,\bar{f}) \in V_1$, we have
$$t(f) - s(f) = d\bar{f}.$$

Next we define $T$ on morphisms.  Suppose $\phi \maps C
\rightarrow C'$ is a chain map between 2-term chain complexes:
$$\xymatrix{
    C_{1} \ar[rr]^{d}
    \ar[dd]_{\phi_{1}} &&
    C_{0} \ar[dd]^{\phi_{0}} \\ \\
    C_{1}' \ar[rr]^{d'} &&
    C_{0}'
}$$ Let $T(C) = V$ and $T(C')=V'$. Then we define $F = T(\phi)$
where $F \maps V \to V'$ is the linear functor with $F_{0}
=\phi_{0}$ and $F_{1} = \phi_{0} \oplus \phi_{1}$.  To check that
$F$ really is a linear functor, note that it is linear on objects
and morphisms.  Moreover, it preserves the source and target,
identity-assigning and composition maps because all these are
defined in terms of addition and the differential in the chain
complexes $C$ and $C'$, and $\phi$ is linear and preserves the
differential.

We leave it the reader to verify that $T$ and $S$ are indeed
functors.  To show that $S$ and $T$ form an equivalence, we
construct natural isomorphisms $\alpha \maps ST \To 1_{2\Vect}$
and $\beta \maps TS \To 1_{2\Term}$.

To construct $\alpha$, consider a $2$-vector space $V$.
Applying $S$ to $V$ we obtain the $2$-term chain complex
$$\xymatrix{
    \ker(s) \ar[rr]^<<<<<<<<<{t|_{\ker(s)}}
   && V_{0}.
}$$ Applying $T$ to this result, we obtain a $2$-vector space $V'$
with the space $V_{0}$ of objects and the space $V_{0} \oplus
\ker(s)$ of morphisms.  The source map for this 2-vector space is
given by $s'(x,\bar{f}) = x$, the target map is given by
$t'(x,\bar{f}) = x + t(\bar{f})$, and the identity-assigning map
is given by $i'(x) = (x,0)$.  We thus can define an isomorphism
$\alpha_V \maps V' \to V$ by setting
\begin{eqnarray*}
(\alpha_V)_0(x) &=& x , \\
(\alpha_V)_1(x,\bar{f}) &=& i(x) + \bar{f} .
\end{eqnarray*}
It is easy to check that $\alpha_V$ is a linear functor.
It is an isomorphism thanks to
the fact, shown in the proof of Lemma \ref{watereddown},
that every morphism in $V$ can be uniquely written
as $i(x) + \bar{f}$ where $x$ is an object and $\bar{f} \in
\ker(s)$.

To construct $\beta$, consider a $2$-term chain complex, $C$, given by
$$\xymatrix{
   C_{1} \ar[rr]^{d}
    && C_{0}.
}$$
Then $T(C)$ is the $2$-vector space with the space $C_{0}$ of objects,
the space $C_{0} \oplus C_{1}$ of morphisms, together with the
source and target maps $s \maps (x,\bar{f}) \mapsto x$,
$t \maps (x,\bar{f}) \mapsto x + d\bar{f}$ and the identity-assigning
map $i \maps x \mapsto (x,0)$.
Applying the functor $S$ to this $2$-vector space
we obtain a $2$-term chain complex $C'$ given by:
$$\xymatrix{
   \ker(s) \ar[rr]^<<<<<<<<<<<{t|_{\ker(s)}}
    && C_{0}.
}$$
Since $\ker(s) = \{(x,\bar{f}) | x = 0 \} \subseteq C_0 \oplus C_1$,
there is an obvious isomorphism
$ker(s) \cong C_{1}$.  Using this we obtain an isomorphism
$\beta_C \maps C' \to C$ given by:
$$\xymatrix{
    \ker(s) \ar[rr]^{t|_{\ker(s)}}
    \ar[dd]_{\sim} &&
    C_{0} \ar[dd]^{1} \\ \\
    C_{1} \ar[rr]^{d} &&
    C_{0}
}$$
where the square commutes because of how we have defined $t$.

We leave it to the reader to verify that $\alpha$ and $\beta$ are
indeed natural isomorphisms. \qed

As mentioned in the Introduction, the idea behind Proposition
\ref{1-1vs} goes back at least to Grothendieck \cite{Gr}, who
showed that groupoids in the category of abelian groups are
equivalent to 2-term chain complexes of abelian groups. There are
many elaborations of this idea, some of which we will mention
later, but for now the only one we really {\it need} involves
making $2\Vect$ and $2\Term$ into $2$-categories and showing that
they are $2$-equivalent as $2$-categories.  To do this, we require
the notion of a `linear natural transformation' between linear
functors.  This will correspond to a chain homotopy between chain
maps.

\begin{defn} \et Given two linear functors $F,G \maps V \to W$
between 2-vector spaces,
a {\bf linear natural transformation} $\alpha \maps F \To G$
is a natural transformation in $\Vect$.
\end{defn}

\begin{defn} \et \label{22Vect} We define {\bf 2Vect} to be
$\Vect\Cat$, or in other words,
the 2-category of 2-vector spaces, linear functors
and linear natural transformations.
\end{defn}

Recall that in general, given two chain maps $\phi, \psi \maps C
\to C'$, a {\bf chain homotopy} $\tau \maps \phi \To \psi$ is a
family of linear maps $\tau \maps C_{p} \rightarrow C_{p+1}'$ such
that $\tau_{p} d_{p+1}' + d_{p} \tau_{p-1} = \psi_{p} - \phi_{p}$
for all $p$. In the case of $2$-term chain complexes, a chain
homotopy amounts to a map $\tau \maps C_{0} \rightarrow C_{1}'$
satisfying $\tau d' = \psi_{0} - \phi_{0}$ and $d \tau = \psi_{1}
- \phi_{1}$.

\begin{defn} \et \label{22Term}
We define {\bf 2Term} to be the 2-category of
2-term chain complexes, chain maps, and chain homotopies.
\end{defn}

\noindent
We will continue to sometimes use $2\Term$ and $2\Vect$ to stand
for the underlying categories of these (strict) 2-categories.
It will be clear by context whether we mean the category or
the $2$-category.

The next result strengthens Proposition \ref{1-1vs}.

\begin{theorem} \et \label{equivof2vs} The $2$-category {\rm $2$Vect}
is $2$-equivalent to the $2$-category {\rm $2$Term}.
\end{theorem}

\noindent{\bf Proof. }
We begin by constructing 2-functors
$$S\maps {\rm 2Vect} \rightarrow {\rm 2Term}$$
and
$$T\maps {\rm 2Term} \rightarrow {\rm 2Vect}.$$
By Proposition \ref{1-1vs}, we
need only to define $S$ and $T$ on 2-morphisms.
Let $V$ and $V'$ be $2$-vector spaces, $F,G \maps V \to V'$
linear functors, and $\theta \maps F \To G$ a linear natural
transformation.  Then we define the chain homotopy
$S(\theta) \maps S(F) \To S(G)$ by
$$  S(\theta)(x) = \bar{\theta}_x ,$$
using the fact that a 0-chain $x$ of $S(V)$ is the same
as an object $x$ of $V$.
Conversely, let $C$ and $C'$ be 2-term chain complexes,
$\phi,\psi \maps C \to C'$ chain maps and $\tau \maps \phi
\To \psi$ a chain homotopy.  Then we define the linear
natural transformation $T(\tau) \maps T(\phi) \To T(\psi)$
by
$$  T(\tau)(x) = (\phi_{0}(x),\tau(x)), $$
where we use the description of a morphism in $S(C')$
as a pair consisting of its source and its arrow part,
which is a 1-chain in $C'$.  We leave it to the reader to check
that $S$ is really a chain homotopy, $T$ is really a linear natural
transformation, and that the natural isomorphisms
$\alpha \maps ST \To 1_{2\Vect}$
and $\beta \maps TS \To 1_{2\Term}$ defined in the
proof of Proposition \ref{1-1vs} extend to this 2-categorical
context.
\qed

We conclude this section with a little categorified linear algebra.
We consider the direct sum and tensor product of 2-vector spaces.

\begin{proposition} \et Given $2$-vector spaces $V=(V_{0}, V_{1},
s,t,i ,\circ)$ and $V'=(V_{0}', V_{1}',$ $s',t',i',\circ')$, there is
a $2$-vector space $V \oplus V'$ having:
\begin{itemize}
\item $V_0 \oplus V'_0$ as its vector space of objects,
\item $V_1 \oplus V'_1$ as its vector space of morphisms,
\item $s\oplus s'$ as its source map,
\item $t \oplus t'$ as its target map,
\item $i \oplus i'$ as its identity-assigning map, and
\item $\circ \oplus \circ'$ as its composition map.
\end{itemize}
\end{proposition}

\noindent{\bf Proof. }  The proof amounts to a routine
verification that the diagrams in Definition \ref{co} commute.
\qed

\begin{proposition} \label{tenprod}
\et Given $2$-vector spaces $V=(V_{0}, V_{1}, s,t,i ,\circ)$
and $V'=(V_{0}', V_{1}',$ $s',t',i',\circ')$, there is a
$2$-vector space $V \otimes V'$ having:
\begin{itemize}
\item $V_0 \otimes V'_0$ as its vector space of objects,
\item $V_1 \otimes V'_1$ as its vector space of morphisms,
\item $s\otimes s'$ as its source map,
\item $t \otimes t'$ as its target map,
\item $i \otimes i'$ as its identity-assigning map, and
\item $\circ \otimes \circ'$ as its composition map.
\end{itemize}
\end{proposition}

\noindent{\bf Proof. }  Again, the proof is a routine verification.
\qed

We now check the correctness of the above definitions by showing
the universal properties of the direct sum and tensor product.
These universal properties only require the category structure of
$2\Vect$, not its 2-category structure, since the necessary
diagrams commute `on the nose' rather than merely up to a
2-isomorphism, and uniqueness holds up to isomorphism, not just up
to equivalence. The direct sum is what category theorists call a
`biproduct': both a product and coproduct, in a compatible way
\cite{Mac}:

\begin{proposition} \et The direct sum $V \oplus V'$ is
the biproduct of the 2-vector spaces $V$ and $V'$, with
the obvious inclusions
\[    i \maps V \to V \oplus V', \qquad i' \maps V' \to V \oplus V' \]
and projections
\[   p \maps V \oplus V' \to V, \qquad p' \maps V \oplus V' \to V' .\]
\end{proposition}

\noindent{\bf Proof. } A routine verification. \qed

Since the direct sum $V \oplus V'$ is a product in the
categorical sense, we may also denote it by $V \times V'$,
as we do now in defining a `bilinear functor',
which is used in stating the universal property of the tensor
product:

\begin{defn} \et \label{bilinearfunct} Let $V, V',$ and $W$ be
$2$-vector spaces.  A {\bf bilinear functor}
$F \maps V \times V' \rightarrow W$
is a functor such that the underlying map on objects
$$F_{0}\maps V_0 \times V'_0 \rightarrow W_{0}$$ and
the underlying map on morphisms
$$F_{1}\maps V_1 \times V'_1 \rightarrow W_{1}$$
are bilinear.
\end{defn}

\begin{proposition} \et Let $V, V',$ and $W$ be $2$-vector spaces.  Given
a bilinear functor $F\maps V \times V' \rightarrow W$ there exists a unique
linear functor $\tilde{F}\maps V \otimes V' \rightarrow W$ such that
$$\xymatrix{
   V \times V'
   \ar[rr]^{F}
   \ar[dd]_{i}
   &&
   W
   \\ \\
   V \otimes V'
   \ar[rruu]_{\tilde{F}}
}$$
commutes, where $i \maps V \times V' \rightarrow V \otimes V'$ is
given by $(v,w) \mapsto v \otimes w$
for $(v,w) \in (V \times V')_{0}$ and
$(f,g) \mapsto f \otimes g$
for $(f,g) \in (V \times V')_{1}$.
\end{proposition}

\noindent{\bf Proof. } The existence and uniqueness of
$\tilde{F_0} \maps (V \otimes V')_0 \to W_0$
and \break
$\tilde{F_1} \maps (V \otimes V')_1 \to W_1$
follow from the universal property of the tensor product
of vector spaces, and it is then straightforward to check
that $\tilde{F}$ is a linear functor.  \qed

We can also form the tensor product of linear functors.
Given linear functors $F \maps V \to V'$ and $G \maps W \to W'$,
we define $F \otimes G \maps V \otimes V' \to W \otimes W'$
by setting:
\[
\begin{array}{ccl}
(F \otimes G)_{0} &=& F_{0} \otimes G_{0}, \\
(F \otimes G)_{1} &=& F_{1} \otimes G_{1}.
\end{array}
\]
Furthermore, there is an `identity object' for the tensor
product of 2-vector spaces.    In $\Vect$, the ground
field $k$ acts as the identity for tensor product:
there are canonical isomorphisms $k \otimes V \cong V$
and $V \otimes k \cong V$.  For 2-vector spaces, a
categorified version of the ground field plays this role:

\begin{proposition} \et There exists a unique $2$-vector space $K$,
the {\bf categorified ground field}, with
$K_{0} = K_{1} = k$ and $s,t,i = 1_{k}.$
\end{proposition}

\noindent{\bf Proof. }  Lemma \ref{watereddown}
implies that there is a unique way to define composition in $K$
making it into a 2-vector space.  In fact, every morphism
in $K$ is an identity morphism.  \qed

\begin{proposition} \et Given any $2$-vector space $V$,
there is an isomorphism $\ell_V \maps K \otimes V \to V$,
which is defined on objects by $a \tensor v \mapsto av$ and on
morphisms by $a \tensor f \mapsto af$.   There is also an
isomorphism $r_V \maps V \otimes K \to V$, defined similarly.
\end{proposition}

\noindent{\bf Proof. }  This is straightforward.
\qed

The functors $\ell_V$ and $r_V$ are a categorified version
of left and right multiplication by scalars.  Our 2-vector
spaces also have a categorified version of addition, namely
a linear functor
\[     + \maps V \oplus V \to V \]
mapping any pair $(x,y)$ of objects or morphisms to $x+y$.
Combining this with scalar multiplication by the object
$-1 \in K$, we obtain another linear functor
\[      - \maps V \oplus V \to V \]
mapping $(x,y)$ to $x-y$.  This is the sense in which our 2-vector
spaces are equipped with a categorified version of subtraction.
All the usual rules governing addition of vectors, subtraction of
vectors, and scalar multiplication hold `on the nose' as
equations.

One can show that with the above tensor product, the
category $2\Vect$ becomes a symmetric monoidal category.
One can go further and make the 2-category version of $2\Vect$
into a symmetric monoidal 2-category \cite{DS}, but we will not need
this here.  Now that we have a definition of 2-vector space
and some basic tools of categorified linear algebra we
may proceed to the main focus of this paper:  the definition of a
categorified Lie algebra.


\section{Semistrict Lie $2$-algebras} \label{Lie2algs}
\subsection{Definitions} \label{definitions}

We now introduce the concept of a `Lie $2$-algebra', which blends
together the notion of a Lie algebra with that of a category.  As
mentioned previously, to obtain a Lie 2-algebra
we begin with a $2$-vector space and equip
it with a bracket {\it functor}, which satisfies the Jacobi
identity {\it up to a natural isomorphism}, the `Jacobiator'.
Then we require that the Jacobiator satisfy a new coherence
law of its own, the `Jacobiator identity'.  We shall assume
the bracket is bilinear in the sense of Definition \ref{bilinearfunct},
and also skew-symmetric:

\begin{defn} \et Let $V$ and $W$ be 2-vector spaces.
A bilinear functor $F\maps V \times V \to W$ is {\bf
skew-symmetric} if $F(x,y) = -F(y,x)$ whenever $(x,y)$ is an
object or morphism of $V \times V$.  If this is the case we
also say the corresponding linear functor $\tilde{F} \maps V
\otimes V \to W$ is skew-symmetric.
\end{defn}

\noindent
We shall also assume that the Jacobiator is trilinear
and completely antisymmetric:

\begin{defn} \et Let $V$ and $W$ be 2-vector spaces.
A functor $F \maps V^n \to W$ is \textbf{n}-{\bf linear} if
$F(x_1,\dots, x_n)$ is linear in each argument, where $(x_1,
\dots, x_n)$ is an object or morphism of $V^n$.  Given $n$-linear
functors $F,G \maps V^n \to W$, a natural transformation $\theta
\maps F \To G$ is \textbf{n}-{\bf linear} if
$\theta_{x_1,\dots,x_n}$ depends linearly on each object $x_i$,
and {\bf completely antisymmetric} if the arrow part of
$\theta_{x_1,\dots,x_n}$ is completely antisymmetric under
permutations of the objects.
\end{defn}

\noindent
Since we do not weaken the bilinearity or skew-symmetry of the
bracket, we call the resulting sort of Lie 2-algebra `semistrict':

\begin{defn} \et \label{defnlie2alg}
A {\bf semistrict Lie $2$-algebra} consists of:

\begin{itemize}
\item a $2$-vector space $L$
\end{itemize}
equipped with
\begin{itemize}
\item a skew-symmetric
bilinear functor, the {\bf bracket}, $[\cdot, \cdot]\maps L \times L
\rightarrow L$
\item a completely antisymmetric
trilinear natural isomorphism, the \break {\bf Jacobiator},
$$J_{x,y,z} \maps [[x,y],z] \to [x,[y,z]] + [[x,z],y],$$
\end{itemize}
that is required to satisfy
\begin{itemize}
 \item the {\bf Jacobiator identity}:
    $$J_{[w,x],y,z} ([J_{w,x,z},y] + 1) (J_{w, [x,z], y} +
       J_{[w,z],x,y} + J_{w,x, [y,z]}) = $$
$$[J_{w,x,y},z] (J_{[w,y],x,z} + J_{w, [x,y],z}) ([J_{w,y,z},x] + 1)
([w, J_{x,y,z}] + 1)$$
\end{itemize}
for all $w,x,y,z \in L_{0}$.  (There is only one choice of
identity morphism that can be added to each term to make the
composite well-defined.)
\end{defn}

The Jacobiator identity looks quite intimidating at first. But if
we draw it as a commutative diagram, we see that it relates two
ways of using the Jacobiator to rebracket the expression
$[[[w,x],y],z]$:

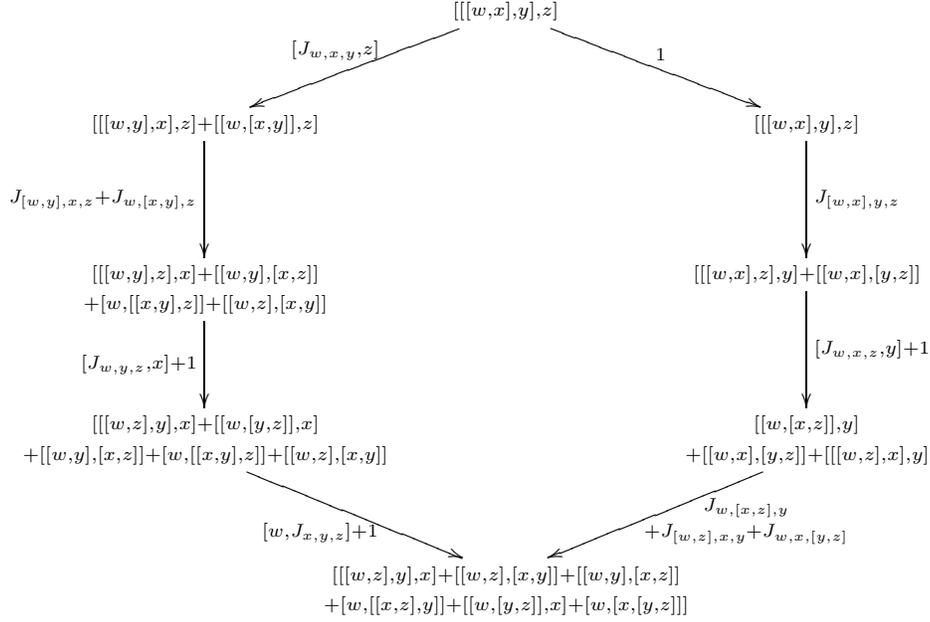
$$ \def\objectstyle{\scriptstyle}
  \def\labelstyle{\scriptstyle}
   \xy
   (0,35)*+{[[[w,x],y],z]}="1";
   (-40,20)*+{[[[w,y],x],z] + [[w,[x,y]],z]}="2";
   (40,20)*+{[[[w,x],y],z]}="3";
   (-40,0)*+{[[[w,y],z],x] + [[w,y],[x,z]]}="4'";
   (-40,-4)*+{+ [w,[[x,y],z]] + [[w,z],[x,y]]}="4";
   (40,0)*+{[[[w,x],z],y] + [[w,x],[y,z]]}="5'";
   (-40,-20)*+{[[[w,z],y],x] + [[w,[y,z]],x]}="6'";
   (-40,-24)*+{+ [[w,y], [x,z]] + [w,[[x,y],z]] + [[w,z],[x,y]]}="6";
   (40,-20)*+{[[w,[x,z]],y]}="7'";
   (40,-24)*+{+ [[w,x],[y,z]] + [[[w,z],x],y]}="7";
   (0,-40)*+{[[[w,z],y],x] + [[w,z],[x,y]]  + [[w,y],[x,z]]}="8'";
   (0,-44)*+{+ [w,[[x,z],y]]  + [[w,[y,z]],x] + [w,[x,[y,z]]]}="8";
            (32,-31)*{J_{w,[x,z],y} }; 
            (32,-34.5)*{+  J_{[w,z],x,y} + J_{w,x,[y,z]}};
        {\ar_{[J_{w,x,y},z]}                   "1";"2"};
        {\ar^{1}                               "1";"3"};
        {\ar_{J_{[w,y],x,z} + J_{w,[x,y],z}}   "2";"4'"};
        {\ar_{[J_{w,y,z},x]+1}                 "4";"6'"};
        {\ar^{J_{[w,x],y,z}}                   "3";"5'"};
        {\ar^{[J_{w,x,z},y]+1}                 "5'";"7'"};
        {\ar_{[w,J_{x,y,z}]+1 \; \; }          "6";"8'"};
        {\ar^{}                                "7";"8'"};
\endxy
\\ \\
$$

\vskip 1em

\noindent
Here the identity morphisms come from terms on which we are
not performing any manipulation.
The reader will surely be puzzled by the fact that we
have included an identity morphism along one edge of this
commutative octagon.  This is explained in the next section, where
we show that the Jacobiator identity is really just a
disguised version of the `Zamolodchikov tetrahedron
equation', which plays an important role in the theory of
higher-dimensional knots and braided monoidal 2-categories
\cite{BN,Crans,DS,KV}.  The Zamolochikov tetrahedron
equation says that two 2-morphisms are equal, each of which
is the vertical composite of four factors.  However, when
we translate this equation into the language of Lie 2-algebras,
one of these factors is an identity 2-morphism.

In the rest of this paper, the term `Lie 2-algebra' will always
refer to a semistrict one as defined above.   We continue by setting
up a 2-category of these Lie 2-algebras.  A homomorphism between
Lie $2$-algebras should preserve both the $2$-vector space
structure and the bracket. However, we shall require that it
preserve the bracket only {\it up to isomorphism}
--- or more precisely, up to a natural isomorphism satisfying a
suitable coherence law.  Thus, we make the following definition.

\begin{defn} \et \label{lie2algfunct} Given
Lie $2$-algebras $L$ and $L'$, a {\bf homomorphism} $F \maps L
\rightarrow L'$ consists of:

\begin{itemize}
  \item A linear functor $F$ from the underlying $2$-vector space of $L$
    to that of $L'$, and

  \item a skew-symmetric bilinear natural transformation
    $$F_{2}(x,y)\maps [F_{0}(x), F_{0}(y)] \rightarrow F_{0}[x,y]$$
\end{itemize}

such that the following diagram commutes:

$$\xymatrix{
     [[F_0(x), F_0(y)], F_0(z)]
        \ar[rrrr]^<<<<<<<<<<<<<<<<<<<<<<{J_{F_0(x), F_0(y),
F_0(z)}}        \ar[dd]_{[F_{2}, 1]}
         &&&& [F_0(x), [F_0(y), F_0(z)]] + [[F_0(x),F_0(z)], F_0(y)]
        \ar[dd]^{[1, F_2] + [F_2, 1]} \\ \\
         [F_0[x,y], F_0(z)]
        \ar[dd]_{F_{2}}
         &&&& [F_0(x), F_0[y,z]] + [F_0[x,z], F_0(y)]
        \ar[dd]^{F_{2} + F_{2}} \\ \\
         F_0[[x,y],z]
        \ar[rrrr]^{F_1(J_{x,y,z})}
         &&&& F_0[x,[y,z]] + F_0[[x,z],y]}$$
         \end{defn}

\noindent
Here and elsewhere we omit the arguments of natural
transformations such as $F_2$ and $G_2$ when these are obvious from
context.

We also have `2-homomorphisms' between homomorphisms:

\begin{defn} \et \label{lie2algnattrans} Let $F,G \maps
L \to L'$ be Lie 2-algebra homomorphisms.  A {\bf 2-homomorphism}
$\theta \maps F \To G$ is a linear natural transformation from $F$
to $G$ such that the following diagram commutes:

$$\xymatrix{
    [F_{0}(x), F_{0}(y)]
     \ar[rr]^{F_{2}}
     \ar[dd]_{[\theta_{x}, \theta_{y}]}
      && F_{0}[x,y]
     \ar[dd]^{\theta_{[x,y]}} \\ \\
      [G_{0}(x), G_{0}(y)]
     \ar[rr]^{G_{2}}
      && G_{0}[x,y] }$$

\end{defn}

\noindent Definitions \ref{lie2algfunct} and \ref{lie2algnattrans}
are closely modelled after the usual definitions of `monoidal
functor' and `monoidal natural transformation' \cite{Mac}.

Next we introduce composition and
identities for homomorphisms and 2-homomorphisms.
The composite of a pair of Lie $2$-algebra homomorphisms
$F\maps L \rightarrow L'$ and $G\maps L' \rightarrow L''$ is
given by letting the functor $FG \maps L \to L''$ be the usual
composite of $F$ and $G$:

$$\xymatrix{
    L
     \ar[rr]^{F}
      &&  L'
     \ar[rr]^{G}
      && L''}$$

\noindent while letting $(FG)_{2}$ be defined as the following
composite:

$$\xymatrix{
    [(FG)_{0}(x), (FG)_{0}(y)]
     \ar[rr]^<<<<<<<<<{(FG)_{2}}
     \ar[dd]_{G_{2}}
      && (FG)_{0}[x,y]  \\ \\
      G_{0}[F_{0}(x), F_{0}(y)]
     \ar[uurr]_{F_{2}\circ G}}.
$$

\noindent
where $F_2 \circ G$ is the result of whiskering the
functor $G$ by the natural transformation $F_2$.
The identity homomorphism $1_L \maps L \to L$ has
the identity functor as its underlying functor, together with
an identity natural transformation as $(1_L)_2$.
Since 2-homomorphisms are just natural transformations
with an extra property, we vertically and horizontally compose
these the usual way, and an identity 2-homomorphism is just
an identity natural transformation.  We obtain:

\begin{proposition} \et There is a strict 2-category {\bf Lie2Alg}
with semistrict Lie $2$-algebras as objects, homomorphisms
between these as morphisms, and $2$-homo-\break morphisms between those
as 2-morphisms, with composition and identities defined as above.
\end{proposition}

\noindent{\bf Proof. } We leave it to the reader to check the
details, including that the composite of homomorphisms is
a homomorphism, this composition is associative,
and the vertical and horizontal composites of
2-homomorphisms are again 2-homomorphisms.  \qed

Finally, note that there is a forgetful $2$-functor from {\rm
Lie$2$Alg} to {\rm $2$Vect}, which is analogous to the
forgetful functor from {\rm LieAlg} to {\rm Vect}.

\subsection{Relation to Topology} \label{topology}

The key novel feature of a Lie 2-algebra is the coherence law for
the Jacobiator: the so-called `Jacobiator identity' in Definition
\ref{defnlie2alg}. At first glance this identity seems rather
arcane. In this section, we `explain' this identity by showing its
relation to the Zamolodchikov tetrahedron equation. This equation
plays a role in the theory of knotted surfaces in 4-space which is
closely analogous to that played by the Yang--Baxter equation, or
third Reidemeister move, in the theory of ordinary knots in
3-space.  In fact, we shall see that just as any Lie algebra gives
a solution of the Yang--Baxter equation, any Lie 2-algebra gives a
solution of the Zamolodchikov tetrahedron equation.

We begin by recalling the Yang--Baxter equation:

\begin{defn} \et Given a vector space $V$ and an
isomorphism $B \maps V \otimes V \to V \otimes V$,
we say $B$ is a {\bf Yang--Baxter operator} if
it satisfies the {\bf Yang--Baxter equation},
which says that:
$$(B \otimes 1)(1 \otimes B)(B \otimes 1)
= (1 \otimes B)(B \otimes 1)(1 \otimes B),$$
or in other words, that this diagram commutes:
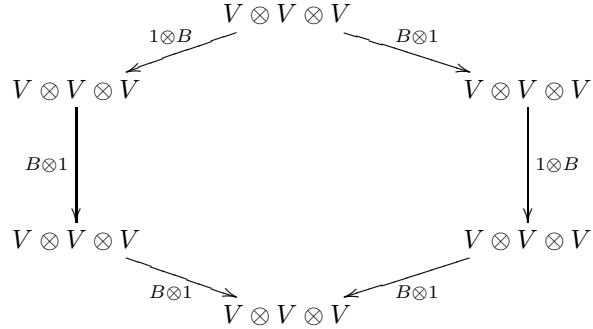
\[
\begin{xy}
    (-2, 20)*+{V \otimes V \otimes V}="1";
    (-30,10)*+{V \otimes V \otimes V}="2";
    (-30,-10)*+{V \otimes V \otimes V}="3";
    (-2,-20)*+{V \otimes V \otimes V}="4";
    (30,-10)*+{V \otimes V \otimes V}="5";
    (30, 10)*+{V \otimes V \otimes V}="6";
        {\ar^{B \otimes 1} "1";"6"};
        {\ar_{1 \otimes B} "1";"2"};
        {\ar_{B \otimes 1} "2";"3"};
        {\ar_{B \otimes 1} "3";"4"};
        {\ar^{B \otimes 1} "5";"4"};
        {\ar^{1 \otimes B} "6";"5"};
\end{xy}
\]
\end{defn}

If we draw $B \maps V \otimes V \to V \otimes V $ as a braiding:
$$
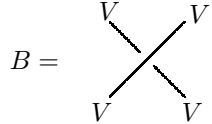
\begin{xy} \xoverv~{(-5,5)}{(5,5)}{(-5,-5)}{(5,-5)};
    (-5,6.4)*{V};
    (7,6)*{V};
    (-15,0)*{B=};
    (-6,-7)*{V};
    (6,-7)*{V};
\end{xy}$$
the Yang--Baxter equation says that:
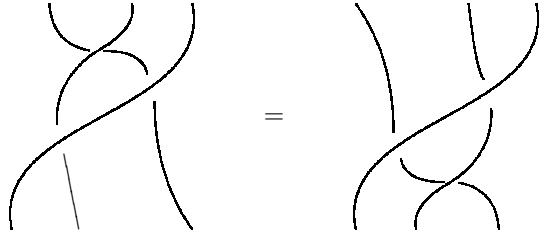
\[
\def\objectstyle{\scriptstyle}
\def\labelstyle{\scriptstyle}
  \xy
   (12,15)*{}="C";
   (4,15)*{}="B";
   (-7,15)*{}="A";
   (12,-15)*{}="3";
   (-3,-15)*{}="2";
   (-12,-15)*{}="1";
       "C";"1" **\crv{(15,0)& (-15,0)};
       (-5,-5)*{}="2'";
       (7,2)*{}="3'";
     "2'";"2" **\crv{};
     "3'";"3" **\crv{(7,-8)};
       \vtwist~{"A"}{"B"}{(-6,-1)}{(6,5.7)};
\endxy
 \qquad = \qquad
   \xy
   (-12,-15)*{}="C";
   (-4,-15)*{}="B";
   (7,-15)*{}="A";
   (-12,15)*{}="3";
   (3,15)*{}="2";
   (12,15)*{}="1";
       "C";"1" **\crv{(-15,0)& (15,0)};
       (5,5)*{}="2'";
       (-7,-2)*{}="3'";
     "2'";"2" **\crv{(4,6)};
     "3'";"3" **\crv{(-7,8)};
       \vtwist~{"A"}{"B"}{(6,1)}{(-6,-5.7)};
\endxy
\]

\noindent This is called the `third Reidemeister move' in knot
theory \cite{BZ}, and it gives the most important
relations in Artin's presentation of the braid group
\cite{Birman}. As a result, any solution of the Yang--Baxter
equation gives an invariant of braids.

In general, almost any process of switching the order of
two things can be thought of as a `braiding'.  This idea
is formalized in the concept of a braided monoidal category,
where the braiding is an isomorphism
\[    B_{x,y} \maps x \tensor y \to y \tensor x .\]
Since the bracket $[x,y]$ in a Lie algebra measures the
difference between $xy$ and $yx$, it should not be too surprising
that we can get a Yang--Baxter operator from any Lie algebra.
And since the third Reidemeister move involves three strands,
while the Jacobi identity involves three Lie algebra elements,
it should also not be surprising that the Yang--Baxter equation is
actually {\it equivalent} to the Jacobi identity in a suitable
context:

\begin{proposition} \et Let $L$ be a vector space equipped
with a skew-symmetric bilinear operation $[\cdot,\cdot] \maps L
\times L \to L$.  Let $L' = k \oplus L$ and define the isomorphism
$B \maps L' \tensor L' \to L' \tensor L'$ by
$$B((a,x) \otimes (b,y))
= (b,y) \otimes (a,x) + (1,0) \otimes (0, [x,y]).$$
Then $B$ is a solution of the Yang--Baxter equation if
and only if $[\cdot,\cdot]$ satisfies the Jacobi identity.
\end{proposition}

\noindent {\bf Proof. }  The proof is a calculation
best left to the reader.  \qed

The nice thing is that this result has a
higher-dimensional analogue, obtained by categorifying
everything in sight!  The analogue of the Yang--Baxter
equation is called the `Zamolodchikov tetrahedron equation':

\begin{defn} \et Given a $2$-vector space $V$ and an
invertible linear functor
$B \maps V \otimes V \to V \otimes V$, a linear natural
isomorphism
$$Y \maps (B \otimes 1)(1 \otimes B)(B \otimes 1) \To
          (1 \otimes B)(B \otimes 1)(1 \otimes B) $$
satisfies the {\bf Zamolodchikov tetrahedron equation} if
$$
[
(Y \otimes 1)
\circ
(1 \otimes 1 \otimes B)
(1 \otimes B \otimes 1)
(B \otimes 1 \otimes 1)
] 
[
(1 \otimes B \otimes 1)
(B \otimes 1 \otimes 1)
\circ 
(1 \otimes Y) 
\circ 
(B \otimes 1 \otimes 1) 
]
$$
$$
[
(1 \otimes B \otimes 1)
(1 \otimes 1 \otimes B)
\circ 
(Y \otimes 1) 
\circ 
(1 \otimes 1 \otimes B)
]
[
(1 \otimes Y)
\circ 
(B \otimes 1 \otimes 1)
(1 \otimes B \otimes 1)
(1 \otimes 1 \otimes B)
]
$$
$$=$$
$$
[
(B \otimes 1 \otimes 1)
(1 \otimes B \otimes 1)
(1 \otimes 1 \otimes B)
\circ 
(Y \otimes 1) 
]
[
(B \otimes 1 \otimes 1)
\circ 
(1 \otimes Y)
\circ 
(B \otimes 1 \otimes 1) 
(1 \otimes B \otimes 1)
]
$$
$$
[
(1 \otimes 1 \otimes B)
\circ 
(Y \otimes 1)
\circ 
(1 \otimes 1 \otimes B) 
(1 \otimes B \otimes 1)
]
[
(1 \otimes 1 \otimes B)
(1 \otimes B \otimes 1)
(B \otimes 1 \otimes 1)
\circ 
(1 \otimes Y) 
]
,$$
where $\circ$ represents the whiskering of a linear functor by
a linear natural transformation.
\end{defn}

To see the significance of this complex but beautifully
symmetrical equation, one should think of
$Y$ as the surface in 4-space traced out by
the {\it process of performing} the third Reidemeister
move:

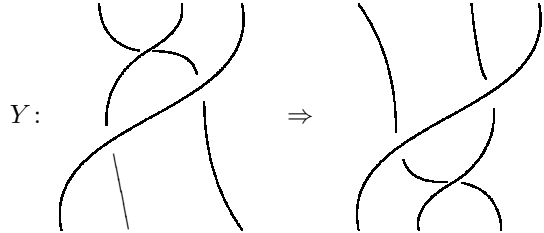
\[   Y \maps
\def\objectstyle{\scriptstyle}
\def\labelstyle{\scriptstyle}
  \xy
   (12,15)*{}="C";
   (4,15)*{}="B";
   (-7,15)*{}="A";
   (12,-15)*{}="3";
   (-3,-15)*{}="2";
   (-12,-15)*{}="1";
       "C";"1" **\crv{(15,0)& (-15,0)};
       (-5,-5)*{}="2'";
       (7,2)*{}="3'";
     "2'";"2" **\crv{};
     "3'";"3" **\crv{(7,-8)};
       \vtwist~{"A"}{"B"}{(-6,-1)}{(6,5.7)};
\endxy
 \quad \To \quad \xy
   (-12,-15)*{}="C";
   (-4,-15)*{}="B";
   (7,-15)*{}="A";
   (-12,15)*{}="3";
   (3,15)*{}="2";
   (12,15)*{}="1";
       "C";"1" **\crv{(-15,0)& (15,0)};
       (5,5)*{}="2'";
       (-7,-2)*{}="3'";
     "2'";"2" **\crv{(4,6)};
     "3'";"3" **\crv{(-7,8)};
       \vtwist~{"A"}{"B"}{(6,1)}{(-6,-5.7)};
\endxy
\]

\noindent
Then the Zamolodchikov tetrahedron equation says the surface
traced out by first performing
the third Reidemeister move on a threefold crossing and then
sliding the result under a fourth strand:
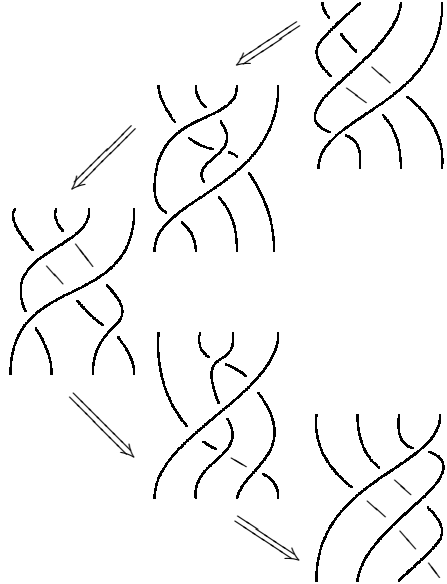
\[
  \xy 0;/r.13pc/:   
    (0,50)*+{       
 \xy 
 (-15,-20)*{}="T1";
 (-5,-20)*{}="T2";
 (5,-20)*{}="T3";
 (15,-20)*{}="T4";
 (-14,20)*{}="B1";
 (-5,20)*{}="B2";
 (5,20)*{}="B3";
 (15,20)*{}="B4";
    "T1"; "B4" **\crv{(-15,-7) & (15,-5)}
        \POS?(.25)*{\hole}="2x" \POS?(.47)*{\hole}="2y" \POS?(.6)*{\hole}="2z";
    "T2";"2x" **\crv{(-4,-12)};
    "T3";"2y" **\crv{(5,-10)};
    "T4";"2z" **\crv{(16,-9)};
 (-15,-5)*{}="3x";
    "2x"; "3x" **\crv{(-18,-10)};
    "3x"; "B3" **\crv{(-13,0) & (4,10)}
        \POS?(.3)*{\hole}="4x" \POS?(.53)*{\hole}="4y";
    "2y"; "4x" **\crv{};
    "2z"; "4y" **\crv{};
 (-15,10)*{}="5x";
    "4x";"5x" **\crv{(-17,6)};
    "5x";"B2" **\crv{(-14,12)}
        \POS?(.6)*{\hole}="6x";
    "6x";"B1" **\crv{(-14,18)};
    "4y";"6x" **\crv{(-8,10)};
 \endxy
    }="T";
    (-40,30)*+{
 \xy 
 (-15,-20)*{}="b1";
 (-5,-20)*{}="b2";
 (5,-20)*{}="b3";
 (14,-20)*{}="b4";
 (-14,20)*{}="T1";
 (-5,20)*{}="T2";
 (5,20)*{}="T3";
 (15,20)*{}="T4";
    "b1"; "T4" **\crv{(-15,-7) & (15,-5)}
        \POS?(.25)*{\hole}="2x" \POS?(.47)*{\hole}="2y" \POS?(.65)*{\hole}="2z";
    "b2";"2x" **\crv{(-5,-15)};
    "b3";"2y" **\crv{(5,-10)};
    "b4";"2z" **\crv{(14,-9)};
 (-15,-5)*{}="3x";
    "2x"; "3x" **\crv{(-15,-10)};
    "3x"; "T3" **\crv{(-15,15) & (5,10)}
        \POS?(.38)*{\hole}="4y" \POS?(.65)*{\hole}="4z";
    "T1";"4y" **\crv{(-14,16)};
    "T2";"4z" **\crv{(-5,16)};
    "2y";"4z" **\crv{(-10,3) & (10,2)} \POS?(.6)*{\hole}="5z";
    "4y";"5z" **\crv{(-5,5)};
    "5z";"2z" **\crv{(5,4)};
 \endxy
    }="TL";
    (-75,0)*+{
 \xy 
 (-14,20)*{}="T1";
 (-4,20)*{}="T2";
 (4,20)*{}="T3";
 (15,20)*{}="T4";
 (-15,-20)*{}="B1";
 (-5,-20)*{}="B2";
 (5,-20)*{}="B3";
 (15,-20)*{}="B4";
    "B1";"T4" **\crv{(-15,5) & (15,-5)}
        \POS?(.25)*{\hole}="2x" \POS?(.49)*{\hole}="2y" \POS?(.65)*{\hole}="2z";
    "2x";"T3" **\crv{(-20,10) & (5,10) }
        \POS?(.45)*{\hole}="3y" \POS?(.7)*{\hole}="3z";
    "2x";"B2" **\crv{(-5,-14)};
        "T1";"3y" **\crv{(-16,17)};
        "T2";"3z" **\crv{(-5,17)};
        "3z";"2z" **\crv{};
        "3y";"2y" **\crv{};
        "B3";"2z" **\crv{ (5,-5) &(20,-10)}
            \POS?(.4)*{\hole}="4z";
        "2y";"4z" **\crv{(6,-8)};
        "4z";"B4" **\crv{(15,-15)};
 \endxy
    }="ML";
    (-40,-30)*+{
 \xy 
 (-14,20)*{}="T1";
 (-4,20)*{}="T2";
 (4,20)*{}="T3";
 (15,20)*{}="T4";
 (-15,-20)*{}="B1";
 (-5,-20)*{}="B2";
 (5,-20)*{}="B3";
 (15,-20)*{}="B4";
    "B1";"T4" **\crv{(-15,-5) & (15,5)}
        \POS?(.38)*{\hole}="2x" \POS?(.53)*{\hole}="2y" \POS?(.7)*{\hole}="2z";
    "T1";"2x" **\crv{(-15,5)};
    "2y";"B2" **\crv{(10,-10) & (-6,-10)}
        \POS?(.45)*{\hole}="4x";
    "2z";"B3" **\crv{ (15,0)&(15,-10) & (6,-16)}
        \POS?(.7)*{\hole}="5x";
    "T3";"2y" **\crv{(5,10)& (-6,18) }
        \POS?(.5)*{\hole}="3x";
    "T2";"3x" **\crv{(-5,15)};
    "3x";"2z" **\crv{(7,11)};
    "2x";"4x" **\crv{(-3,-7)};
    "4x";"5x" **\crv{};
    "5x";"B4" **\crv{(15,-15)};
 \endxy
    }="BL";
(0,-50)*+{
 \xy 
 (15,20)*{}="T1";
 (5,20)*{}="T2";
 (-5,20)*{}="T3";
 (-15,20)*{}="T4";
 (15,-20)*{}="B1";
 (5,-20)*{}="B2";
 (-5,-20)*{}="B3";
 (-15,-20)*{}="B4";
    "T1"; "B4" **\crv{(15,7) & (-15,5)}
        \POS?(.25)*{\hole}="2x" \POS?(.45)*{\hole}="2y" \POS?(.6)*{\hole}="2z";
    "T2";"2x" **\crv{(4,12)};
    "T3";"2y" **\crv{(-5,10)};
    "T4";"2z" **\crv{(-16,9)};
 (15,5)*{}="3x";
    "2x"; "3x" **\crv{(18,10)};
    "3x"; "B3" **\crv{(13,0) & (-4,-10)}
        \POS?(.3)*{\hole}="4x" \POS?(.53)*{\hole}="4y";
    "2y"; "4x" **\crv{};
    "2z"; "4y" **\crv{};
 (15,-10)*{}="5x";
    "4x";"5x" **\crv{(17,-6)};
    "5x";"B2" **\crv{(14,-12)}
        \POS?(.6)*{\hole}="6x";
    "6x";"B1" **\crv{};
    "4y";"6x" **\crv{};
 \endxy
    }="B";
            (-20,65)*{}="X1";
            (-35,55)*{}="X2";               
                {\ar@{=>} "X1";"X2"};       
            (-60,40)*{}="X1";               
            (-75,25)*{}="X2";               
                {\ar@{=>} "X1";"X2"};
            (-60,-40)*{}="X2";
            (-75,-25)*{}="X1";
                {\ar@{=>} "X1";"X2"};
            (-20,-65)*{}="X2";
            (-35,-55)*{}="X1";
                {\ar@{=>} "X1";"X2"};
  \endxy
\]
is isotopic to that traced out by first sliding the threefold crossing
under the fourth strand and then performing the third
Reidemeister move:
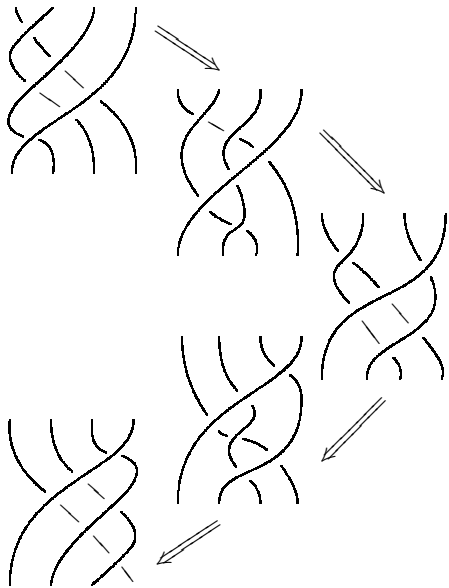
\[
  \xy 0;/r.13pc/:   
    (0,50)*+{       
 \xy 
 (-15,-20)*{}="T1";
 (-5,-20)*{}="T2";
 (5,-20)*{}="T3";
 (15,-20)*{}="T4";
 (-14,20)*{}="B1";
 (-5,20)*{}="B2";
 (5,20)*{}="B3";
 (15,20)*{}="B4";
    "T1"; "B4" **\crv{(-15,-7) & (15,-5)}
        \POS?(.25)*{\hole}="2x" \POS?(.47)*{\hole}="2y" \POS?(.6)*{\hole}="2z";
    "T2";"2x" **\crv{(-4,-12)};
    "T3";"2y" **\crv{(5,-10)};
    "T4";"2z" **\crv{(16,-9)};
 (-15,-5)*{}="3x";
    "2x"; "3x" **\crv{(-18,-10)};
    "3x"; "B3" **\crv{(-13,0) & (4,10)}
        \POS?(.3)*{\hole}="4x" \POS?(.53)*{\hole}="4y";
    "2y"; "4x" **\crv{};
    "2z"; "4y" **\crv{};
 (-15,10)*{}="5x";
    "4x";"5x" **\crv{(-17,6)};
    "5x";"B2" **\crv{(-14,12)}
        \POS?(.6)*{\hole}="6x";
    "6x";"B1" **\crv{(-14,18)};
    "4y";"6x" **\crv{(-8,10)};
 \endxy
    }="T";
(40,30)*+{
 \xy 
 (14,-20)*{}="T1";
 (4,-20)*{}="T2";
 (-4,-20)*{}="T3";
 (-15,-20)*{}="T4";
 (15,20)*{}="B1";
 (5,20)*{}="B2";
 (-5,20)*{}="B3";
 (-15,20)*{}="B4";
    "B1";"T4" **\crv{(15,5) & (-15,-5)}
        \POS?(.38)*{\hole}="2x" \POS?(.53)*{\hole}="2y" \POS?(.7)*{\hole}="2z";
    "T1";"2x" **\crv{(15,-5)};
    "2y";"B2" **\crv{(-10,10) & (6,10)}
        \POS?(.45)*{\hole}="4x";
    "2z";"B3" **\crv{ (-15,0)&(-15,10) & (-6,16)}
        \POS?(.7)*{\hole}="5x";
    "T3";"2y" **\crv{(-5,-10)& (6,-18) }
        \POS?(.5)*{\hole}="3x";
    "T2";"3x" **\crv{(5,-15)};
    "3x";"2z" **\crv{(-7,-11)};
    "2x";"4x" **\crv{(3,7)};
    "4x";"5x" **\crv{};
    "5x";"B4" **\crv{(-15,15)};
 \endxy
    }="TR";
    (75,0)*+{
 \xy 
 (14,-20)*{}="T1";
 (4,-20)*{}="T2";
 (-4,-20)*{}="T3";
 (-15,-20)*{}="T4";
 (15,20)*{}="B1";
 (5,20)*{}="B2";
 (-5,20)*{}="B3";
 (-15,20)*{}="B4";
    "B1";"T4" **\crv{(15,-5) & (-15,5)}
        \POS?(.25)*{\hole}="2x" \POS?(.49)*{\hole}="2y" \POS?(.65)*{\hole}="2z";
    "2x";"T3" **\crv{(20,-10) & (-5,-10) }
        \POS?(.45)*{\hole}="3y" \POS?(.7)*{\hole}="3z";
    "2x";"B2" **\crv{(5,14)};
        "T1";"3y" **\crv{(16,-17)};
        "T2";"3z" **\crv{(5,-17)};
        "3z";"2z" **\crv{};
        "3y";"2y" **\crv{};
        "B3";"2z" **\crv{ (-5,5) &(-20,10)}
            \POS?(.4)*{\hole}="4z";
        "2y";"4z" **\crv{(-6,8)};
        "4z";"B4" **\crv{(-15,15)};
 \endxy
    }="MR";
    (40,-30)*+{
 \xy 
 (15,20)*{}="b1";
 (5,20)*{}="b2";
 (-5,20)*{}="b3";
 (-14,20)*{}="b4";
 (14,-20)*{}="T1";
 (5,-20)*{}="T2";
 (-5,-20)*{}="T3";
 (-15,-20)*{}="T4";
    "b1"; "T4" **\crv{(15,7) & (-15,5)}
        \POS?(.25)*{\hole}="2x" \POS?(.47)*{\hole}="2y" \POS?(.65)*{\hole}="2z";
    "b2";"2x" **\crv{(5,15)};
    "b3";"2y" **\crv{(-5,10)};
    "b4";"2z" **\crv{(-14,9)};
 (15,5)*{}="3x";
    "2x"; "3x" **\crv{(15,10)};
    "3x"; "T3" **\crv{(15,-15) & (-5,-10)}
        \POS?(.38)*{\hole}="4y" \POS?(.65)*{\hole}="4z";
    "T1";"4y" **\crv{(14,-16)};
    "T2";"4z" **\crv{(5,-16)};
    "2y";"4z" **\crv{(10,-3) & (-10,-2)} \POS?(.6)*{\hole}="5z";
    "4y";"5z" **\crv{(5,-5)};
    "5z";"2z" **\crv{(-5,-4)};
 \endxy
    }="BR";
    (0,-50)*+{
 \xy 
 (15,20)*{}="T1";
 (5,20)*{}="T2";
 (-5,20)*{}="T3";
 (-15,20)*{}="T4";
 (15,-20)*{}="B1";
 (5,-20)*{}="B2";
 (-5,-20)*{}="B3";
 (-15,-20)*{}="B4";
    "T1"; "B4" **\crv{(15,7) & (-15,5)}
        \POS?(.25)*{\hole}="2x" \POS?(.45)*{\hole}="2y" \POS?(.6)*{\hole}="2z";
    "T2";"2x" **\crv{(4,12)};
    "T3";"2y" **\crv{(-5,10)};
    "T4";"2z" **\crv{(-16,9)};
 (15,5)*{}="3x";
    "2x"; "3x" **\crv{(18,10)};
    "3x"; "B3" **\crv{(13,0) & (-4,-10)}
        \POS?(.3)*{\hole}="4x" \POS?(.53)*{\hole}="4y";
    "2y"; "4x" **\crv{};
    "2z"; "4y" **\crv{};
 (15,-10)*{}="5x";
    "4x";"5x" **\crv{(17,-6)};
    "5x";"B2" **\crv{(14,-12)}
        \POS?(.6)*{\hole}="6x";
    "6x";"B1" **\crv{};
    "4y";"6x" **\crv{};
 \endxy
    }="B";
            (20,65)*{}="X1";                
            (35,55)*{}="X2";                
                {\ar@{=>} "X1";"X2"};       
            (60,40)*{}="X1";                
            (75,25)*{}="X2";                
                {\ar@{=>} "X1";"X2"};       
            (60,-40)*{}="X2";
            (75,-25)*{}="X1";
                {\ar@{=>} "X1";"X2"};
            (20,-65)*{}="X2";
            (35,-55)*{}="X1";
                {\ar@{=>} "X1";"X2"};
  \endxy
\]
In short, the Zamolodchikov tetrahedron equation is a
formalization of this commutative octagon:
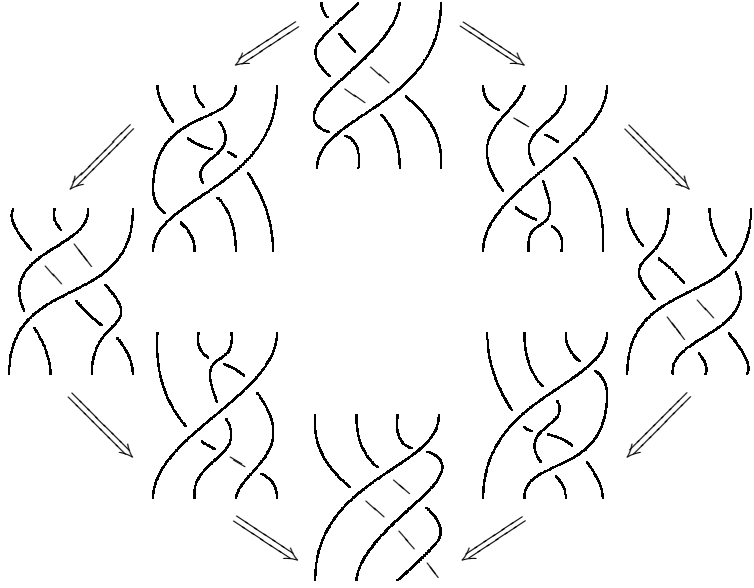
\[
  \xy 0;/r.13pc/:   
    (0,50)*+{       
 \xy 
 (-15,-20)*{}="T1";
 (-5,-20)*{}="T2";
 (5,-20)*{}="T3";
 (15,-20)*{}="T4";
 (-14,20)*{}="B1";
 (-5,20)*{}="B2";
 (5,20)*{}="B3";
 (15,20)*{}="B4";
    "T1"; "B4" **\crv{(-15,-7) & (15,-5)}
        \POS?(.25)*{\hole}="2x" \POS?(.47)*{\hole}="2y" \POS?(.6)*{\hole}="2z";
    "T2";"2x" **\crv{(-4,-12)};
    "T3";"2y" **\crv{(5,-10)};
    "T4";"2z" **\crv{(16,-9)};
 (-15,-5)*{}="3x";
    "2x"; "3x" **\crv{(-18,-10)};
    "3x"; "B3" **\crv{(-13,0) & (4,10)}
        \POS?(.3)*{\hole}="4x" \POS?(.53)*{\hole}="4y";
    "2y"; "4x" **\crv{};
    "2z"; "4y" **\crv{};
 (-15,10)*{}="5x";
    "4x";"5x" **\crv{(-17,6)};
    "5x";"B2" **\crv{(-14,12)}
        \POS?(.6)*{\hole}="6x";
    "6x";"B1" **\crv{(-14,18)};
    "4y";"6x" **\crv{(-8,10)};
 \endxy
    }="T";
    (-40,30)*+{
 \xy 
 (-15,-20)*{}="b1";
 (-5,-20)*{}="b2";
 (5,-20)*{}="b3";
 (14,-20)*{}="b4";
 (-14,20)*{}="T1";
 (-5,20)*{}="T2";
 (5,20)*{}="T3";
 (15,20)*{}="T4";
    "b1"; "T4" **\crv{(-15,-7) & (15,-5)}
        \POS?(.25)*{\hole}="2x" \POS?(.47)*{\hole}="2y" \POS?(.65)*{\hole}="2z";
    "b2";"2x" **\crv{(-5,-15)};
    "b3";"2y" **\crv{(5,-10)};
    "b4";"2z" **\crv{(14,-9)};
 (-15,-5)*{}="3x";
    "2x"; "3x" **\crv{(-15,-10)};
    "3x"; "T3" **\crv{(-15,15) & (5,10)}
        \POS?(.38)*{\hole}="4y" \POS?(.65)*{\hole}="4z";
    "T1";"4y" **\crv{(-14,16)};
    "T2";"4z" **\crv{(-5,16)};
    "2y";"4z" **\crv{(-10,3) & (10,2)} \POS?(.6)*{\hole}="5z";
    "4y";"5z" **\crv{(-5,5)};
    "5z";"2z" **\crv{(5,4)};
 \endxy
    }="TL";
    (-75,0)*+{
 \xy 
 (-14,20)*{}="T1";
 (-4,20)*{}="T2";
 (4,20)*{}="T3";
 (15,20)*{}="T4";
 (-15,-20)*{}="B1";
 (-5,-20)*{}="B2";
 (5,-20)*{}="B3";
 (15,-20)*{}="B4";
    "B1";"T4" **\crv{(-15,5) & (15,-5)}
        \POS?(.25)*{\hole}="2x" \POS?(.49)*{\hole}="2y" \POS?(.65)*{\hole}="2z";
    "2x";"T3" **\crv{(-20,10) & (5,10) }
        \POS?(.45)*{\hole}="3y" \POS?(.7)*{\hole}="3z";
    "2x";"B2" **\crv{(-5,-14)};
        "T1";"3y" **\crv{(-16,17)};
        "T2";"3z" **\crv{(-5,17)};
        "3z";"2z" **\crv{};
        "3y";"2y" **\crv{};
        "B3";"2z" **\crv{ (5,-5) &(20,-10)}
            \POS?(.4)*{\hole}="4z";
        "2y";"4z" **\crv{(6,-8)};
        "4z";"B4" **\crv{(15,-15)};
 \endxy
    }="ML";
    (-40,-30)*+{
 \xy 
 (-14,20)*{}="T1";
 (-4,20)*{}="T2";
 (4,20)*{}="T3";
 (15,20)*{}="T4";
 (-15,-20)*{}="B1";
 (-5,-20)*{}="B2";
 (5,-20)*{}="B3";
 (15,-20)*{}="B4";
    "B1";"T4" **\crv{(-15,-5) & (15,5)}
        \POS?(.38)*{\hole}="2x" \POS?(.53)*{\hole}="2y" \POS?(.7)*{\hole}="2z";
    "T1";"2x" **\crv{(-15,5)};
    "2y";"B2" **\crv{(10,-10) & (-6,-10)}
        \POS?(.45)*{\hole}="4x";
    "2z";"B3" **\crv{ (15,0)&(15,-10) & (6,-16)}
        \POS?(.7)*{\hole}="5x";
    "T3";"2y" **\crv{(5,10)& (-6,18) }
        \POS?(.5)*{\hole}="3x";
    "T2";"3x" **\crv{(-5,15)};
    "3x";"2z" **\crv{(7,11)};
    "2x";"4x" **\crv{(-3,-7)};
    "4x";"5x" **\crv{};
    "5x";"B4" **\crv{(15,-15)};
 \endxy
    }="BL";
    (40,30)*+{
 \xy 
 (14,-20)*{}="T1";
 (4,-20)*{}="T2";
 (-4,-20)*{}="T3";
 (-15,-20)*{}="T4";
 (15,20)*{}="B1";
 (5,20)*{}="B2";
 (-5,20)*{}="B3";
 (-15,20)*{}="B4";
    "B1";"T4" **\crv{(15,5) & (-15,-5)}
        \POS?(.38)*{\hole}="2x" \POS?(.53)*{\hole}="2y" \POS?(.7)*{\hole}="2z";
    "T1";"2x" **\crv{(15,-5)};
    "2y";"B2" **\crv{(-10,10) & (6,10)}
        \POS?(.45)*{\hole}="4x";
    "2z";"B3" **\crv{ (-15,0)&(-15,10) & (-6,16)}
        \POS?(.7)*{\hole}="5x";
    "T3";"2y" **\crv{(-5,-10)& (6,-18) }
        \POS?(.5)*{\hole}="3x";
    "T2";"3x" **\crv{(5,-15)};
    "3x";"2z" **\crv{(-7,-11)};
    "2x";"4x" **\crv{(3,7)};
    "4x";"5x" **\crv{};
    "5x";"B4" **\crv{(-15,15)};
 \endxy
    }="TR";
    (75,0)*+{
 \xy 
 (14,-20)*{}="T1";
 (4,-20)*{}="T2";
 (-4,-20)*{}="T3";
 (-15,-20)*{}="T4";
 (15,20)*{}="B1";
 (5,20)*{}="B2";
 (-5,20)*{}="B3";
 (-15,20)*{}="B4";
    "B1";"T4" **\crv{(15,-5) & (-15,5)}
        \POS?(.25)*{\hole}="2x" \POS?(.49)*{\hole}="2y" \POS?(.65)*{\hole}="2z";
    "2x";"T3" **\crv{(20,-10) & (-5,-10) }
        \POS?(.45)*{\hole}="3y" \POS?(.7)*{\hole}="3z";
    "2x";"B2" **\crv{(5,14)};
        "T1";"3y" **\crv{(16,-17)};
        "T2";"3z" **\crv{(5,-17)};
        "3z";"2z" **\crv{};
        "3y";"2y" **\crv{};
        "B3";"2z" **\crv{ (-5,5) &(-20,10)}
            \POS?(.4)*{\hole}="4z";
        "2y";"4z" **\crv{(-6,8)};
        "4z";"B4" **\crv{(-15,15)};
 \endxy
    }="MR";
    (40,-30)*+{
 \xy 
 (15,20)*{}="b1";
 (5,20)*{}="b2";
 (-5,20)*{}="b3";
 (-14,20)*{}="b4";
 (14,-20)*{}="T1";
 (5,-20)*{}="T2";
 (-5,-20)*{}="T3";
 (-15,-20)*{}="T4";
    "b1"; "T4" **\crv{(15,7) & (-15,5)}
        \POS?(.25)*{\hole}="2x" \POS?(.47)*{\hole}="2y" \POS?(.65)*{\hole}="2z";
    "b2";"2x" **\crv{(5,15)};
    "b3";"2y" **\crv{(-5,10)};
    "b4";"2z" **\crv{(-14,9)};
 (15,5)*{}="3x";
    "2x"; "3x" **\crv{(15,10)};
    "3x"; "T3" **\crv{(15,-15) & (-5,-10)}
        \POS?(.38)*{\hole}="4y" \POS?(.65)*{\hole}="4z";
    "T1";"4y" **\crv{(14,-16)};
    "T2";"4z" **\crv{(5,-16)};
    "2y";"4z" **\crv{(10,-3) & (-10,-2)} \POS?(.6)*{\hole}="5z";
    "4y";"5z" **\crv{(5,-5)};
    "5z";"2z" **\crv{(-5,-4)};
 \endxy
    }="BR";
    (0,-50)*+{
 \xy 
 (15,20)*{}="T1";
 (5,20)*{}="T2";
 (-5,20)*{}="T3";
 (-15,20)*{}="T4";
 (15,-20)*{}="B1";
 (5,-20)*{}="B2";
 (-5,-20)*{}="B3";
 (-15,-20)*{}="B4";
    "T1"; "B4" **\crv{(15,7) & (-15,5)}
        \POS?(.25)*{\hole}="2x" \POS?(.45)*{\hole}="2y" \POS?(.6)*{\hole}="2z";
    "T2";"2x" **\crv{(4,12)};
    "T3";"2y" **\crv{(-5,10)};
    "T4";"2z" **\crv{(-16,9)};
 (15,5)*{}="3x";
    "2x"; "3x" **\crv{(18,10)};
    "3x"; "B3" **\crv{(13,0) & (-4,-10)}
        \POS?(.3)*{\hole}="4x" \POS?(.53)*{\hole}="4y";
    "2y"; "4x" **\crv{};
    "2z"; "4y" **\crv{};
 (15,-10)*{}="5x";
    "4x";"5x" **\crv{(17,-6)};
    "5x";"B2" **\crv{(14,-12)}
        \POS?(.6)*{\hole}="6x";
    "6x";"B1" **\crv{};
    "4y";"6x" **\crv{};
 \endxy
    }="B";
            (-20,65)*{}="X1";
            (-35,55)*{}="X2";               
                {\ar@{=>} "X1";"X2"};       
            (20,65)*{}="X1";                
            (35,55)*{}="X2";                
                {\ar@{=>} "X1";"X2"};       
            (60,40)*{}="X1";                
            (75,25)*{}="X2";                
                {\ar@{=>} "X1";"X2"};       
            (-60,40)*{}="X1";               
            (-75,25)*{}="X2";               
                {\ar@{=>} "X1";"X2"};
            (-60,-40)*{}="X2";
            (-75,-25)*{}="X1";
                {\ar@{=>} "X1";"X2"};
            (60,-40)*{}="X2";
            (75,-25)*{}="X1";
                {\ar@{=>} "X1";"X2"};
            (-20,-65)*{}="X2";
            (-35,-55)*{}="X1";
                {\ar@{=>} "X1";"X2"};
            (20,-65)*{}="X2";
            (35,-55)*{}="X1";
                {\ar@{=>} "X1";"X2"};
  \endxy
\]

\noindent in a 2-category whose 2-morphisms are isotopies
of surfaces in 4-space --- or more precisely, `2-braids'.
Details can be found in HDA1, HDA4 and a number of other
references, going back to the work of Kapranov and Voevodsky
\cite{BLan,BN,CS,Crans,KV}.

In Section \ref{definitions}, we drew the Jacobiator identity
as a commutative octagon.  In fact, that commutative octagon
becomes {\it equivalent} to the octagon for the
Zamolodchikov tetrahedron equation in the following context:

\begin{theorem} \et \label{Zameqn} Let $L$ be a $2$-vector space,
let $[\cdot, \cdot] \maps L \times L \to L$ be a skew-symmetric
bilinear functor, and let $J$ be a completely antisymmetric trilinear
natural transformation with $J_{x,y,z}\maps [[x,y],z]
\to [x,[y,z]] + [[x,z],y]$.  Let $L' = K \oplus L$,
where $K$ is the categorified ground field. Let $B \maps L' \otimes
L' \to L' \otimes L'$ be defined as follows:
$$B((a,x) \otimes (b,y))
= (b,y) \otimes (a,x) + (1,0) \otimes (0, [x,y]).$$
whenever $(a,x)$ and $(b,y)$ are both either objects or morphisms
in $L'$.
Finally, let
$$Y \maps (B \otimes 1)(1 \otimes B)(B \otimes 1) \To
          (1 \otimes B)(B \otimes 1)(1 \otimes B) $$
be defined as follows:
$$Y = (p \otimes p \otimes p) \circ J \circ j$$
where $p \maps L' \to L$ is the projection functor given by the
fact that $L' = K \oplus L$ and
$$j \maps L \rightarrow L' \otimes L' \otimes L'$$
is the linear functor defined by
$$j(x) = (1,0) \otimes (1,0) \otimes (0,x),$$
where $x$ is either an object or morphism of $L$. Then $Y$ is a
solution of the Zamolodchikov tetrahedron equation if and only if
$J$ satisfies the Jacobiator identity.
\end{theorem}

\noindent{\bf Proof. }  Equivalently, we must show that $Y$
satisfies the Zamolodchikov tetrahedron equation if and only if
$J$ satisfies the Jacobiator identity.  Applying the
left-hand side of the Zamolodchikov tetrahedron equation
to an object $(a,w) \otimes (b,x) \otimes (c,y) \otimes (d,z)$
of $L' \otimes L' \otimes L' \otimes L'$ yields
an expression consisting of various uninteresting terms together
with one involving
$$J_{[w,x],y,z} ([J_{w,x,z},y]+1) (J_{w, [x,z], y} +
       J_{[w,z],x,y} + J_{w,x, [y,z]}), $$
while applying the right-hand side produces an expression
with the same uninteresting terms, but also one involving
$$[J_{w,x,y},z] (J_{[w,y],x,z} + J_{w, [x,y],z}) ([J_{w,y,z},x]+1)
([w, J_{x,y,z}]+1) $$
in precisely the same way.  Thus, the two sides are equal if and
only if the Jacobiator identity holds.  The detailed
calculation is quite lengthy. \qed

\begin{corollary} \et  If $L$ is a Lie 2-algebra,
then $Y$ defined as in Theorem \ref{Zameqn} is a solution of the
Zamolodchikov tetrahedron equation.
\end{corollary}

We continue by exhibiting the correlation between our semistrict
Lie $2$-algebras and special versions of Stasheff's
$L_{\infty}$-algebras.

\subsection{$L_{\infty}$-algebras} \label{Linftyalgs}

An $L_\infty$-algebra is a chain complex equipped with a
bilinear skew-symmetric bracket operation that satisfies
the Jacobi identity `up to coherent homotopy'.  In other
words, this identity holds up to a specified chain homotopy,
which in turn satisfies its own identity up to a specified
chain homotopy, and so on {\it ad infinitum}.  Such structures
are also called are `strongly homotopy Lie algebras' or
`sh Lie algebras' for short.  Though their precursors
existed in the literature beforehand, they made their first
notable appearance in a 1985 paper on deformation theory by
Schlessinger and Stasheff \cite{SS}. Since then, they have been
systematically explored and applied in a number of other
contexts \cite{LM, LS, Mar, P}.

Since 2-vector spaces are equivalent to 2-term chain complexes,
as described in Section \ref{2vs}, it should not be surprising
that $L_\infty$-algebras are related to the categorified
Lie algebras we are discussing here.  An elegant but
rather highbrow way to approach this is to use the theory of
operads \cite{MSS}.  An $L_\infty$-algebra is actually an algebra
of a certain operad in the symmetric monoidal category of chain
complexes, called the `$L_\infty$ operad'.  Just as categories
in $\Vect$ are equivalent to 2-term chain complexes, strict
$\omega$-categories in $\Vect$ can be shown equivalent to general
chain complexes, by a similar argument \cite{BH}.   Using this
equivalence, we can transfer the $L_\infty$ operad from the world
of chain complexes to the world of strict $\omega$-category objects
in $\Vect$, and define a {\bf semistrict Lie $\omega$-algebra} to be
an algebra of the resulting operad.

In more concrete terms, a semistrict Lie $\omega$-algebra
is a strict $\omega$-category $L$ having a vector space $L_j$
of $j$-morphisms for all $j \ge 0$, with all source, target and
composition maps being linear.  Furthermore, it is equipped with a
skew-symmetric bilinear bracket functor
\[      [\cdot,\cdot] \maps L \times L \to L  \]
which satisfies the Jacobi identity up to a completely
antisymmetric trilinear natural
isomorphism, the `Jacobiator', which in turn satisfies
the Jacobiator identity up to a completely antisymmetric
quadrilinear modification...
and so on.  By the equivalence mentioned above, such a thing is
really just {\it another way of looking at} an $L_\infty$-algebra.

Using this, one can show that a semistrict Lie $\omega$-algebra
with only identity $j$-morphisms for $j > 1$ is the same as a semistrict
Lie 2-algebra!  But luckily, we can prove a result along these lines
without using or even mentioning the concepts of `operad',
`$\omega$-category' and the like.  Instead, for the sake of an
accessible presentation, we shall simply recall the definition
of an $L_\infty$-algebra and prove that the 2-category of
semistrict Lie 2-algebras is equivalent to a 2-category of
`2-term' $L_\infty$-algebras: that is, those having a
zero-dimensional space of $j$-chains for $j > 1$.

Henceforth, all algebraic objects mentioned are considered over
a fixed field $k$ of characteristic other than 2.  We make
consistent use of the usual sign convention when dealing
with graded objects.  That
is, whenever we interchange something of degree $p$ with something
of degree $q,$ we introduce a sign of $(-1)^{pq}.$ The following
conventions regarding graded vector spaces, permutations,
unshuffles, etc., follow those of Lada and Markl \cite{LM}.

For graded indeterminates $x_{1}, \ldots, x_{n}$ and a permutation
$\sigma \in S_{n}$ we define the {\bf Koszul sign}
$\epsilon (\sigma) = \epsilon(\sigma; x_{1}, \dots, x_{n})$ by
$$x_{1} \wedge \cdots \wedge x_{n} = \epsilon(\sigma; x_{1},
\ldots, x_{n}) \cdot x_{\sigma(1)} \wedge \cdots \wedge
x_{\sigma(n)},$$ which must be satisfied in the free graded-commutative
algebra on $x_{1}, \ldots, x_{n}.$  This is
nothing more than a formalization of what has already been said
above.  Furthermore, we define
$$\chi(\sigma) = \chi(\sigma;
x_{1}, \dots, x_{n}) := \textrm{sgn} (\sigma) \cdot
\epsilon(\sigma; x_{1}, \dots, x_{n}).$$ Thus, $\chi(\sigma)$
takes into account the sign of the permutation in $S_{n}$ and the
sign obtained from iteration of the basic convention.

If $n$ is a natural number and $1 \leq j \leq n-1$ we say that
$\sigma \in S_{n}$ is an $(j,n-j)${\bf -unshuffle} if
$$ \sigma(1) \leq\sigma(2) \leq \cdots \leq \sigma(j)
\hspace{.2in} \textrm{and} \hspace{.2in} \sigma(j+1) \leq
\sigma(j+2) \leq \cdots \leq \sigma(n).$$ Readers familiar with
shuffles will recognize unshuffles as their inverses. A
\emph{shuffle} of two ordered sets (such as a deck of cards) is a
permutation of the ordered union preserving the order of each of
the given subsets. An \emph{unshuffle} reverses this process.  A
simple example should clear up any confusion:

\begin{example} \et When $n=3$, the $(1,2)$-unshuffles in $S_{3}$ are:
$$\emph{id} =  \begin{pmatrix}
  1 & 2 & 3 \\
  1 & 2 & 3
\end{pmatrix}, \hspace{.2in}
(132) = \begin{pmatrix}
  1 & 2 & 3 \\
  3 & 1 & 2
\end{pmatrix}, \hspace{.1in} \emph{and} \hspace{.1in}
(12) = \begin{pmatrix}
  1 & 2 & 3 \\
  2 & 1 & 3
\end{pmatrix}.$$
\end{example}

The following definition of an $L_{\infty}$-structure was
formulated by Stasheff in 1985, see \cite{SS}.  This definition
will play an important role in what will follow.

\begin{defn} \et \label{L-alg} An
{\bf $\mathbf{L_{\infty}}$-algebra} is a graded vector space $V$
equipped with a system $\{l_{k}| 1 \leq k < \infty\}$ of linear maps
$l_{k} \maps V^{\otimes k} \rightarrow V$ with $\deg(l_{k}) = k-2$
which are totally antisymmetric in the sense that
\begin{eqnarray}
   l_{k}(x_{\sigma(1)}, \dots,x_{\sigma(k)}) =
   \chi(\sigma)l_{k}(x_{1}, \dots, x_{n})
\label{antisymmetry}
\end{eqnarray}
for all $\sigma \in S_{n}$ and $x_{1}, \dots, x_{n} \in V,$ and,
moreover, the following generalized form of the Jacobi identity
holds for $0 \le n < \infty :$
\begin{eqnarray}
   \displaystyle{\sum_{i+j = n+1}
   \sum_{\sigma}
   \chi(\sigma)(-1)^{i(j-1)} l_{j}
   (l_{i}(x_{\sigma(1)}, \dots, x_{\sigma(i)}), x_{\sigma(i+1)},
   \ldots, x_{\sigma(n)}) =0,}
\label{megajacobi}
\end{eqnarray}
where the summation is taken over all $(i,n-i)$-unshuffles with $i
\geq 1.$
\end{defn}

While somewhat puzzling at first, this definition truly does
combine the important aspects of Lie algebras and chain complexes.
The map $l_1$ makes $V$ into a chain complex, since this map has
degree $-1$ and equation (\ref{megajacobi}) says its square is
zero. Moreover, the map $l_{2}$ resembles a Lie bracket, since it
is skew-symmetric in the graded sense by equation
(\ref{antisymmetry}). In what follows, we usually denote $l_1(x)$
as $dx$ and $l_2(x,y)$ as $[x,y]$.  The higher $l_k$ maps are
related to the Jacobiator, the Jacobiator identity, and the higher
coherence laws that would appear upon further categorification of
the Lie algebra concept.

To make this more precise, let us
refer to an $L_{\infty}$-algebra $V$ with $V_{n} = 0$ for
$n \geq k$ as a \textbf{\textit{k}-term $\mathbf{L_{\infty}}$-algebra.}
Note that a $1$-term $L_{\infty}$-algebra is simply an ordinary Lie
algebra, where $l_{3} =0$ gives the
Jacobi identity.  However, in a $2$-term $L_{\infty}$-algebra,
we no longer have a trivial $l_{3}$ map.  Instead, equation
(\ref{megajacobi}) says that the Jacobi identity for the
0-chains $x,y,z$ holds up to a term of the form
$dl_3(x,y,z)$.  We do, however, have $l_{4} = 0$, which
provides us with the coherence law that $l_{3}$ must satisfy.

Since we will be making frequent use of these $2$-term
$L_{\infty}$-algebras, it will be advantageous to keep track of
their ingredients.

\begin{lemma} \et \label{rmk}
A $2$-term $L_{\infty}$-algebra, $V,$ consists of the following
data:

\begin{itemize}
  \item two vector spaces $V_{0}$ and
   $V_{1}$ together with a linear map $d\maps V_{1} \rightarrow V_{0},$

  \item a bilinear map $l_{2}\maps V_{i} \times V_{j}
   \rightarrow V_{i+j},$ where $0 \le i+j \le 1$,
   \newline which we denote more suggestively as $[\cdot, \cdot],$

  \item a trilinear map $l_{3}\maps V_{0} \times V_{0} \times
   V_{0} \rightarrow V_{1}.$

\end{itemize}

These maps satisfy:

\begin{itemize}
  \item[(a)] $[x,y] = -[y,x]$,
  \item[(b)] $[x,h] = -[h,x]$,
  \item[(c)] $[h,k]=0$,
  \item[(d)] $l_{3}(x,y,z)$ is totally antisymmetric in the
    arguments $x,y,z$,
  \item[(e)] $d([x,h]) = [x,dh]$,
  \item[(f)] $[dh,k] = [h,dk]$,
  \item[(g)] $d(l_{3}(x,y,z)) = -[[x,y],z] + [[x,z],y] + [x,[y,z]]$,
  \item[(h)] $l_{3}(dh,x,y) = - [[x,y],h] + [[x,h],y] + [x,[y,h]]$,
  \item[(i)] $[l_{3}(w,x,y),z]
    + [l_{3}(w,y,z),x] + l_{3}([w,y],x,z) + l_{3}([x,z],w,y) =$
    $$[l_{3}(w,x,z),y] + [l_{3}(x,y,z),w] + l_{3}([w,x],y,z) + $$
    $$ l_{3}([w,z],x,y)+ l_{3}([x,y],w,z)+ l_{3}([y,z],w,x),$$
\end{itemize}
for all $w,x,y,z \in V_{0}$ and $h, k \in V_{1}.$
\end{lemma}

\noindent{\rm Proof. } Note that
$(a)-(d)$ hold by equation (\ref{antisymmetry}) of Definition
\ref{L-alg} while $(e)-(i)$ follow from (\ref{megajacobi}).
\qed

We notice that $(a)$ and $(b)$ are the usual skew-symmetric
properties satisfied by the bracket in a Lie algebra; $(c)$ arises
simply because there are no 2-chains.  Equations $(e)$ and $(f)$
tell us how the differential and bracket interact, while
$(g)$ says that the Jacobi identity no longer holds on the nose,
but up to chain homotopy.  We will use $(g)$ to define the
Jacobiator in the Lie 2-algebra corresponding to a 2-term
$L_\infty$-algebra.  Equation $(h)$ will give the naturality
of the Jacobiator.  Similarly, $(i)$ will give the Jacobiator
identity.

We continue by defining homomorphisms between $2$-term
$L_{\infty}$-algebras:

\begin{defn} \et \label{Linftyhomo}
Let $V$ and $V'$ be $2$-term $L_{\infty}$-algebras. An
\textbf{$\mathbf{L_{\infty}}$-homomorphism}
$\phi \maps V \rightarrow V'$ consists of:

\begin{itemize}
  \item a chain map $\phi \maps V \to V'$ (which consists of
  linear maps $\phi_0 \maps V_0 \to V'_0$ and
              $\phi_1 \maps V_1 \to V'_1$ preserving the differential),
  \item a skew-symmetric bilinear map $\phi_{2} \maps V_{0} \times V_{0} \to
  V_{1}'$,
\end{itemize}

such that the following equations hold for all $x,y,z \in V_0$, $h
\in V_{1}:$

\begin{itemize}
\item $d (\phi_{2}(x,y)) = \phi_{0}[x,y] - [\phi_{0}(x), \phi_{0}(y)]$
\item $\phi_{2}(x,dh) = \phi_{1}[x,h] - [\phi_{0}(x), \phi_{1}(h)]$
\item $ [\phi_2(x,y), \phi_0(z)]  + \phi_2([x,y],z) + \phi_1(l_3(x,y,z))
= $ $ l_3(\phi_0(x),\phi_0(y), \phi_0(z))  + [\phi_0(x), \phi_2(y,z)] + [\phi_2(x,z), \phi_0(y)] + \phi_2(x, [y,z]) + \phi_2([x,z],y)$


\end{itemize}
\end{defn}

\noindent The reader should note the similarity between this
definition and that of homomorphisms between Lie 2-algebras
(Definition \ref{lie2algfunct}).  In particular, the first two
equations say that $\phi_{2}$ defines a chain homotopy from
$[\phi(\cdot), \phi(\cdot)]$ to $\phi[\cdot, \cdot]$, where these
are regarded as chain maps from $V \otimes V$ to $V'$.  The third
equation in the above definition is just a chain complex version
of the commutative square in Definition \ref{lie2algfunct}.

To make 2-term $L_\infty$-algebras and $L_\infty$-homomorphisms
between them into a category, we must describe composition and identities.
We compose a pair of $L_{\infty}$-homomorphisms $\phi \maps V
\rightarrow V'$ and $\psi \maps V' \rightarrow V''$ by letting the
chain map $\phi \psi \maps V \to V''$ be the usual composite:
$$
\xymatrix{
    V
     \ar[rr]^{\phi}
      &&  V'
     \ar[rr]^{\psi}
      && V''}
$$

\noindent while defining $(\phi \psi)_{2}$ as follows:
$$ (\phi \psi)_2(x,y) = \psi_2(\phi_0(x),\phi_0(y)) +
                              \psi_1(\phi_2(x,y)).$$
This is just a chain complex version of how we compose homomorphisms
between Lie 2-algebras.  The identity homomorphism $1_V \maps V \to V$
has the identity chain map as its underlying map, together with
$(1_V)_2 = 0$.

With these definitions, we obtain:

\begin{proposition} \et
There is a category {\bf 2TermL$_\mathbf\infty$} with
2-term $L_{\infty}$-algebras as objects and
$L_\infty$-homomorphisms as morphisms.
\end{proposition}

\noindent{\bf Proof. }  We leave this an exercise for
the reader.    \qed

Next we establish the equivalence between the category of Lie
$2$-algebras and that of $2$-term $L_{\infty}$-algebras.  This result
is based on the equivalence between 2-vector spaces and 2-term chain
complexes described in Proposition \ref{1-1vs}.

\begin{theorem} \et \label{1-1}The categories {\rm Lie$2$Alg}
and {\rm 2TermL$_{\infty}$} are equivalent.
\end{theorem}

\noindent{\bf Proof. }  First we sketch how
to construct a functor $T\maps {\rm 2TermL_{\infty}}
\rightarrow {\rm Lie2Alg}$.  Given a $2$-term $L_{\infty}$-algebra
$V$ we construct the Lie $2$-algebra $L=T(V)$ as follows.

We construct the underlying 2-vector space of $L$ as in
the proof of Proposition \ref{1-1vs}.  Thus $L$ has vector
spaces of objects and morphisms
\begin{eqnarray*}
  L_{0} & = & V_{0}, \\
  L_{1} & = & V_{0} \oplus V_{1},
\end{eqnarray*}
and we denote a morphism $f\maps x
\rightarrow y$ in $L_{1}$ by $f=(x, \bar{f})$
where $x \in V_{0}$ is the source of $f$ and $\bar{f} \in V_{1}$
is its arrow part.  The source, target, and identity-assigning maps in $L$
are given by
\begin{eqnarray*}
  s(f) &=& s(x, \bar{f}) = x ,\\
  t(f) &=& t(x, \bar{f}) = x + d\bar{f}, \\
  i(x) &=& (x, 0),
\end{eqnarray*}
and we have $t(f) - s(f) = d\bar{f}$.  The composite of two morphisms in
$L_{1}$ is given as in the proof of Lemma
\ref{watereddown}.  That is, given $f = (x, \bar{f})\maps x
\rightarrow y,$ and $g = (y, \bar{g})\maps y \rightarrow z$,
we have
$$f g:= (x, \bar{f} + \bar{g}).$$

We continue by equipping $L=T(V)$ with the additional structure
which makes it a Lie $2$-algebra.  First, we use the degree-zero
chain map $l_2 \maps V \otimes V \to V $ to define the bracket
functor $[\cdot, \cdot]\maps L \times L \rightarrow L.$  For a
pair of objects $x,y \in L_0$ we define $[x,y] = l_2(x,y)$, where
we use the `$l_2$' notation in the $L_\infty$-algebra $V$ to avoid
confusion with the bracket in $L$. The bracket functor is
skew-symmetric and bilinear on objects since $l_{2}$ is.  This is
not sufficient, however.  It remains to define the bracket functor
on pairs of morphisms.

We begin by defining the bracket on pairs of morphisms where
one morphism is an identity.  We do this as follows:
given a morphism $f= (x, \bar{f})\maps x \rightarrow y $ in $L_{1}$
and an object $z \in L_{0}$, we define
$$[1_z,f]:= (l_2(z,x), l_2(z, \bar{f})),$$
$$[f,1_z]:= (l_2(x,z), l_2(\bar{f}, z)). $$
Clearly these morphisms have the desired sources;
we now verify that they also have the desired targets.  Using
the fact that $t(f) = s(f) + d\bar{f}$ for any morphism $f \in
L_{1}$, we see that:
\begin{eqnarray*}
t[1_{z},f] &=& s[1_{z},f] + dl_2(z,\bar{f}) \cr
       &=& l_2(z,x) + l_2(z, d\bar{f}) \; \; \; \; \textrm{by $(e)$ of
       Lemma \ref{rmk}} \cr
       &=& l_2(z,x) + l_2(z, y-x) \cr
       &=& l_2(z, y)
\end{eqnarray*}
as desired, using the bilinearity of $l_2$.  Similarly we have
$t[f, 1_{z}] = l_2(y, z)$.

These definitions together with the desired functoriality
of the bracket force us to define the bracket of
an arbitrary pair of morphisms $f\maps x \rightarrow y$,
$g\maps a \rightarrow b$ as follows:
\begin{eqnarray*}
[f,g] &=& [f 1_y, 1_a g] \cr
       &:=& [f,1_a] \, [1_y,g] \cr
       &=& (l_2(x,a), l_2(\bar{f}, a)) \; (l_2(y,a), l_2(y, \bar{g})) \cr
       &=& (l_2(x,a), l_2(\bar{f}, a) + l_2(y, \bar{g})).
\end{eqnarray*}
On the other hand, they also force us to define it as:
\begin{eqnarray*}
[f,g] &=& [1_x f, g 1_b] \\
      &:=& [1_x,g] \, [f,1_b]   \\
      &=& (l_2(x,a), l_2(x,\bar{g})) \; (l_2(x,b), l_2(\bar{f},b)) \\
      &=& (l_2(x,a), l_2(x,\bar{g}) + l_2(\bar{f},b))  .
\end{eqnarray*}
Luckily these are compatible!  We have
\begin{equation}
  l_2(\bar{f}, a) + l_2(y, \bar{g}) =
    l_2(x,\bar{g}) + l_2(\bar{f},b)
\label{magic.equation}
\end{equation}
because the left-hand side minus the right-hand side equals
$l_2(d\bar{f},\bar{g}) - l_2(\bar{f}, d\bar{g})$, which vanishes
by $(f)$ of Lemma \ref{rmk}.

At this point we relax and use $[\cdot,\cdot]$ to stand both for
the bracket on objects in $L$ and the $L_\infty$-algebra $V$. We
will not, however, relax when it comes to the morphisms in $L$
since $[\cdot, \cdot] \neq l_2(\cdot, \cdot)$ even on morphisms
that are arrow parts, that is, morphisms in ker($s$).  By the
above calculations, the bracket of morphisms $f \maps x \to y$, 
$g \maps a \to b$ in $L$ is given by
\begin{eqnarray*}
   [f,g] &=& ([x,a], l_2(\bar{f}, a) + l_2(y, \bar{g}))   \\
         &=& ([x,a], l_2(x,\bar{g}) + l_2(\bar{f},b))  .
\end{eqnarray*}

The bracket $[\cdot,\cdot] \maps L \times L \to L$ is
clearly bilinear on objects.  Either of the above formulas
shows it is also bilinear on morphisms, since the source,
target and arrow part of a morphism depend linearly on the
morphism, and the bracket in $V$ is bilinear.
The bracket is also skew-symmetric: this is clear for
objects, and can be seen for morphisms if we use {\it both} the
above formulas.

To show that $[\cdot, \cdot]\maps L \times L \rightarrow L$ is a
functor, we must check that it preserves identities and
composition.
We first show that $[1_{x}, 1_{y}] = 1_{[x,y]}$, where $x, y \in
L_{0}$.  For this we use the fact that identity morphisms are
precisely those with vanishing arrow part.  Either formula for
the bracket of morphisms gives
\begin{eqnarray*}
  [1_{x}, 1_{y}] &=& ([x,y],0) \cr
                 &=& 1_{[x,y]}.
\end{eqnarray*}

\noindent
To show that the bracket preserves composition, consider the
morphisms $f=(x,\bar{f}), f'= (y, \bar{f'}), g= (a, \bar{g}),$ and
$g'=(b,\bar{g'})$ in $L_{1},$ where $f\maps x \rightarrow y,$ $f'\maps y
\rightarrow z,$ $g\maps a \rightarrow b,$ and $g'\maps b \rightarrow c.$
We must show
$$[f f', g g'] = [f,g] [f',g'].$$
On the one hand, the definitions give
$$
  [f f', g g'] =
        ([x,a], l_2(\bar{f}, a) +l_2(\bar{f'},a) + l_2(z, \bar{g}) + l_2(z, \bar{g'})),
$$
while on the other, they give
$$
  [f,g] [f',g'] =
       ([x,a], l_2(\bar{f},a) + l_2(y, \bar{g}) + l_2(\bar{f'},b) + l_2(z, \bar{g'}))
$$
Therefore, it suffices to show that
$$l_2(\bar{f'},a) + l_2(z,\bar{g}) = l_2(y,\bar{g}) + l_2(\bar{f'},b).$$
After a relabelling of variables, this is just equation
(\ref{magic.equation}).

Next we define the Jacobiator for $L$ and check its properties.
We set
$$J_{x,y,z} := ([[x,y],z], l_{3}(x,y,z)).  $$
Clearly the source of $J_{x,y,z}$ is $[[x,y],z]$ as desired,
while its target is $[x, [y,z]] + [[x,z],y]$ by condition
$(g)$ of Lemma \ref{rmk}.
To show $J_{x,y,z}$ is natural one can check that
is natural in each argument.  We check naturality in the
third variable, leaving the other two as exercises for the
reader.  Let $f \maps z
\rightarrow z'.$ Then, $J_{x,y,z}$ is natural in $z$ if the following
diagram commutes:
$$\xymatrix{
    [[x,y],z]
     \ar[rrr]^{[[1_x,1_y],f]}
     \ar[dd]_{J_{x,y,z}}
      &&& [[x,y],z']
     \ar[dd]^{J_{x,y,z'}} \\ \\
      [[x,z],y]+ [x,[y,z]]
     \ar[rrr]^{[[1_x,f],1_y]+ [1_x,[1_y,f]]}
      &&& [[x,z'],y]+[x,[y,z']] }$$
\\
\noindent Using the formula for the composition and bracket in $L$
this means that we need
$$([[x,y],z], \bar{J}_{x,y,z'} + l_2([x,y],\bar{f})) =
([[x,y],z], l_2(l_2(x, \bar{f}), y) + (x,l_2(y,\bar{f})) +
\bar{J}_{x,y,z}).$$ Thus, it suffices to show that
$$\bar{J}_{x,y,z'} + l_2([x,y],\bar{f}) =
l_2(l_2(x,\bar{f}),y) + l_2(x,l_2(y,\bar{f})) + \bar{J}_{x,y,z}.$$
But $\bar{J}_{x,y,z}$ has been defined as $l_{3}(x,y,z)$ (and
similarly for $\bar{J}_{x,y,z'}$), so now we are required to show
that:
$$l_{3}(x,y,z') + l_2([x,y],\bar{f}) =
l_{3}(x,y,z) + l_2(l_2(x,\bar{f}),y) + l_2(x, l_2(y,\bar{f})),$$
or in other words,
$$ l_2([x,y],\bar{f}) +
l_{3}(x,y,d\bar{f}) = l_2(l_2(x,\bar{f}),y) + l_2(x,
l_2(y,\bar{f})).$$ This holds by condition $(h)$ in Lemma
\ref{rmk} together with the complete antisymmetry of $l_3$.

The trilinearity and complete antisymmetry of the Jacobiator
follow from the corresponding properties of $l_3$.
Finally, condition $(i)$ in Lemma \ref{rmk} gives the Jacobiator
identity:
 $$J_{[w,x],y,z} ([J_{w,x,z},y]+1) (J_{w, [x,z], y} +
       J_{[w,z],x,y} + J_{w,x, [y,z]}) = $$
$$[J_{w,x,y},z] (J_{[w,y],x,z} + J_{w, [x,y],z}) ([J_{w,y,z},x]+1)
([w, J_{x,y,z}]+1).$$

This completes the construction of a Lie 2-algebra $T(V)$ from any
2-term $L_\infty$-algebra $V$.  Next we construct a Lie 2-algebra
homomorphism $T(\phi) \maps T(V) \to T(V')$ from any
$L_\infty$-homomorphism $\phi \maps V \to V'$ between 2-term
$L_\infty$-algebras.

Let $T(V) = L$ and $T(V') = L'$.  We define the underlying linear
functor of $T(\phi)=F$ as in Proposition \ref{1-1vs}.  To make $F$
into a Lie 2-algebra homomorphism we must equip it with a
skew-symmetric bilinear natural transformation $F_{2}$ satisfying
the conditions in Definition \ref{lie2algfunct}. We do this using
the skew-symmetric bilinear map $\phi_{2} \maps V_{0} \times V_{0}
\to V_{1}'$.  In terms of its source and arrow parts, we let
$$F_{2}(x,y) = ([\phi_{0}(x), \phi_{0}(y)], \phi_{2}(x,y)).$$
Computing the target of $F_{2}(x,y)$ we have:
\begin{eqnarray*}
tF_{2}(x,y) &=& sF_{2}(x,y) + d\bar{F_{2}}(x,y) \\
            &=& [\phi_{0}(x), \phi_{0}(y)] + d\phi_{2}(x,y) \\
            &=& [\phi_{0}(x), \phi_{0}(y)] + \phi_{0}[x,y] - [\phi_{0}(x),
            \phi_{0}(y)]\\
            &=& \phi_{0}[x,y] \\
            &=& F_0[x,y]
\end{eqnarray*}
by the first equation in Definition \ref{Linftyhomo} and the fact
that $F_0 = \phi_0$. Thus we have $F_{2}(x,y) \maps [F_{0}(x),
F_{0}(y)] \to F_{0}[x,y]$. Notice that $F_{2}(x,y)$ is bilinear
and skew-symmetric since $\phi_{2}$ and the bracket are. $F_{2}$
is a natural transformation by Theorem \ref{equivof2vs} and the
fact that $\phi_2$ is a chain homotopy from $[\phi(\cdot),
\phi(\cdot)]$ to $\phi([\cdot,\cdot])$, thought of as chain maps
from $V \tensor V$ to $V'$. Finally, the equation in Definition
\ref{Linftyhomo} gives the commutative diagram in Definition
\ref{lie2algfunct}, since the composition of morphisms corresponds
to addition of their arrow parts.

We leave it to the reader to check that $T$ is indeed
a functor.  Next, we describe how to construct a functor $S\maps {\rm
Lie2Alg} \to {\rm 2TermL_{\infty}}$.

Given a Lie 2-algebra $L$ we construct the 2-term $L_\infty$-algebra
$V = S(L)$ as follows.  We define:
\begin{eqnarray*}
  V_{0} &=& L_{0} \\
  V_{1} &=& \ker(s) \subseteq L_{1}.
\end{eqnarray*}
In addition, we define the maps $l_{k}$ as follows:
\begin{itemize}
\item $l_1h = t(h)$ for $h \in V_1 \subseteq L_1$.
\item $l_{2}(x,y) = [x,y]$ for $x,y \in V_0 = L_0$.
\item $l_{2}(x,h) = -l_{2}(h,x) = [1_x, h]$
for $x \in V_0 = L_0$ and $h \in V_1 \subseteq L_1$.
\item $l_2(h,k) = 0$ for $h,k \in V_1 \subseteq L_1$.
\item $l_{3}(x,y,z) = \bar{J}_{x,y,z}$ for $x,y,z \in V_0 = L_0$.
\end{itemize}
As usual, we abbreviate $l_1$ as $d$ and $l_2$ on zero chains as
$[\cdot,\cdot]$.

With these definitions, conditions $(a)$ and $(b)$ of Lemma \ref{rmk}
follow from the antisymmetry of the bracket functor.  Condition $(c)$ is
automatic.  Condition $(d)$ follows from the complete
antisymmetry of the Jacobiator.

For $h \in V_{1}$ and $x \in V_{0}$, the functoriality of $[\cdot,
\cdot]$ gives
\begin{eqnarray*}
   d(l_2(x,h)) &=& t([1_x,h]) \cr
            &=& [t(1_x), t(h)] \cr
            &=& [x, dh],
\end{eqnarray*}
which is $(e)$ of Lemma \ref{rmk}. To obtain $(f)$,
\ref{rmk}, we let $h\maps 0 \rightarrow x$ and $k\maps 0 \rightarrow y$ be
elements of $V_{1}$.  We then consider the following square in $L
\times L$,
$$\xymatrix{
   &
     0
    \ar[rr]^{h}
     && x \\
     0
    \ar[dd]_{k}
     & (0,0)
    \ar[rr]^{(h,1_0)}
    \ar[dd]_{(1_0,k)}
     && (x,0)
    \ar[dd]^{(1_x,k)} \\ \\
     y
     & (0,y)
    \ar[rr]^{(h,1_y)}
     && (x,y)}$$
which commutes by definition of a product category. Since $[\cdot,
\cdot]$ is a functor, it preserves such commutative squares, so
that
$$\xymatrix{
    [0,0]
     \ar[rr]^{[h,1_0]}
     \ar[dd]_{[1_0,k]}
      && [x,0]
     \ar[dd]^{[1_x,k]} \\ \\
      [0,y]
     \ar[rr]^{[h,1_y]}
      && [x,y]}$$
commutes.  Since $[h,1_0]$ and $[1_0,k]$ are easily seen to be
identity morphisms, this implies $[h,1_y] =[1_x,k]$.  This means
that in $V$ we have $l_2(h,y) = l_2(x,k)$, or, since $y$ is the
target of $k$ and $x$ is the target of $h$, simply
$l_2(h,dk)=l_2(dh,k),$ which is $(f)$ of Lemma \ref{rmk}.

Since $J_{x,y,z} \maps [[x,y],z] \to [x,[y,z]] + [[x,z],y]$,
we have
\begin{eqnarray*}
   d(l_{3}(x,y,z)) &=& t(\bar{J}_{x,y,z}) \cr
                   &=& (t-s)(J_{x,y,z}) \cr
                   &=& [x,[y,z]] + [[x,z],y] - [[x,y],z],
\end{eqnarray*}
which gives $(g)$ of Lemma \ref{rmk}.  The naturality of
$J_{x,y,z}$ implies that for any $f \maps z \to z'$, we must have
$$[[1_x,1_y],f] \, J_{x,y,z'} = J_{x,y,z} \,
([[1_x,f],1_y] +[1_x,[1_y,f]]).$$
This implies that in $V$ we have
$$l_2([x,y],\bar{f}) + l_{3}(x,y, z'-z) = l_2(l_2(x,\bar{f}),y) + l_2(x, l_2(y,\bar{f})),$$
for $x,y \in V_{0}$ and $\bar{f} \in V_{1},$ which is $(h)$ of
Lemma \ref{rmk}.

Finally, the Jacobiator identity dictates that the arrow part
of the Jacobiator, $l_3$, satisfies the following equation:
$$l_2(l_{3}(w,x,y),z) + l_2(l_{3}(w,y,z),x) + l_{3}([w,y],x,z) +
l_{3}([x,z],w,y) = $$
$$l_2(l_3 (w,x,z), y) +
 l_2(l_{3}(x,y,z),w) + l_{3}([w,x],y,z) + $$
$$ l_{3}([w,z],x,y)+ l_{3}([x,y],w,z)+ l_{3}([y,z],w,x).$$
This is $(i)$ of Lemma \ref{rmk}.

This completes the construction of a 2-term $L_\infty$-algebra $S(L)$
from any Lie 2-algebra $L$.  Next we construct an
$L_\infty$-homomorphism $S(F) \maps S(L) \to S(L')$ from any
Lie 2-algebra homomorphism $F \maps L \to L'$.

Let $S(L) = V$ and $S(L') = V'$.  We define the underlying chain
map of $S(F) = \phi$ as in Proposition \ref{1-1vs}. To make $\phi$
into an $L_\infty$-homomorphism we must equip it with a
skew-symmetric bilinear map $\phi_{2} \maps V_{0} \times V_{0} \to
V_{1}'$ satisfying the conditions in Definition \ref{Linftyhomo}.
To do this we set
$$
\phi_{2}(x,y) = \bar{F}_{2}(x,y).
$$
The bilinearity and skew-symmetry of $\phi_{2}$
follow from that of $F_{2}$.  Then, since $\phi_{2}$ is
the arrow part of $F_{2}$,
\begin{eqnarray*}
d\phi_{2}(x,y) &=& (t-s)F_{2}(x,y) \\
&=& F_{0}[x,y] - [F_{0}(x), F_{0}(y)] \\
&=& \phi_{0}[x,y] - [\phi_{0}(x), \phi_{0}(y)],
\end{eqnarray*}
by definition of the chain map $\phi$.  The naturality of $F_{2}$
gives the second equation in Definition \ref{Linftyhomo}. Finally,
since composition of morphisms corresponds to addition of arrow
parts, the diagram in Definition \ref{lie2algfunct} gives:
$$l_2(\phi_2(x,y), \phi_0(z))  + \phi_2([x,y],z) + \phi_1(l_3(x,y,z))
=  l_3(\phi_0(x),\phi_0(y), \phi_0(z))  \; + $$ $$ l_2(\phi_0(x), \phi_2(y,z)) + l_2(\phi_2(x,z), \phi_0(y))+ \phi_2(x, [y,z]) + \phi_2([x,z],y),$$
since $\phi_{0} = F_{0}$, $\phi_{1} = F_{1}$ on
elements of $V_{1}$, and the arrow parts of $J$ and $F_{2}$ are
$l_{3}$ and $\phi_{2}$, respectively.

We leave it to the reader to check that
$S$ is indeed a functor, and to construct
natural isomorphisms $\alpha \maps ST \To 1_{\rm Lie2Alg}$ and
$\beta \maps TS \To 1_{{\rm 2TermL}_\infty}$.
\qed


The above theorem also has a 2-categorical version.
We have defined a 2-category of Lie 2-algebras, but not
yet defined a 2-category of 2-term $L_\infty$-algebras.
For this, we need the following:

\begin{defn} \et \label{Linfty2homo} Let $V$ and $V'$ be
2-term $L_\infty$-algebras and let $\phi, \psi \maps V \to
V'$ be $L_{\infty}$-homomorphisms.  An {\bf
$\mathbf L_{\infty}$-2-homomorphism} $\tau \maps \phi \To \psi$ is a
chain homotopy such that the following equation holds for all $x,y \in
V_0$:
\begin{itemize}
\item $\phi_2(x,y) - \psi_2(x,y) =
[\phi_{0}(x), \tau(y)] + [\tau(x), \psi_{0}(y)] - \tau ([x,y])$
\end{itemize}
\end{defn}

\noindent
We leave it as an exercise for the reader to define
vertical and horizontal composition for these 2-homomorphisms,
to define identity 2-homomorphisms, and to prove the following:

\begin{proposition} \et There is a strict 2-category
{\bf 2TermL$_\mathbf\infty$} with 2-term $L_{\infty}$-algebras
as objects, homomorphisms between these as morphisms, and
2-homo- \hfill \break morphisms between those as $2$-morphisms.
\end{proposition}

With these definitions one can strengthen Theorem \ref{1-1}
as follows:

\begin{theorem} \label{1-1'} \et The 2-categories {\rm Lie$2$Alg}
and {\rm 2TermL$_{\infty}$} are 2-equivalent.
\end{theorem}

The main benefit of the results in this section is that they
provide us with another method to create examples of Lie
$2$-algebras.  Instead of thinking of a Lie $2$-algebra as a
category equipped with extra structure, we may work with a
$2$-term chain complex endowed with the structure described in
Lemma \ref{rmk}.  In the next two sections we investigate the
results of trivializing various aspects of a Lie $2$-algebra, or
equivalently of the corresponding 2-term $L_\infty$-algebra.


\section{Strict Lie $2$-algebras} \label{strictlie2algs}

A `strict' Lie $2$-algebra is a categorified version of a Lie
algebra in which all laws hold as equations, not just up to
isomorphism. In a previous paper \cite{B} one of the authors
showed how to construct these starting from `strict Lie 2-groups'.
Here we describe this process in a somewhat more highbrow manner,
and explain how these `strict' notions are special cases of the
semistrict ones described here.

Since we only weakened the Jacobi identity in our definition of
`semistrict' Lie 2-algebra, we need only require that the
Jacobiator be the identity to recover the `strict' notion:

\begin{defn} \et A semistrict Lie 2-algebra is {\bf strict}
if its Jacobiator is the identity.
\end{defn}

\noindent Similarly, requiring that the bracket be strictly
preserved gives the notion of `strict' homomorphism between Lie
2-algebras:

\begin{defn} \et Given semistrict Lie $2$-algebras $L$ and $L'$, a
homomorphism $F \maps L \to L'$ is {\bf strict} if $F_2$ is the
identity.
\end{defn}

\begin{proposition} \et {\rm Lie2Alg}
contains a sub-$2$-category {\bf SLie2Alg} with strict Lie
$2$-algebras as objects, strict homomorphisms between these as
morphisms, and arbitrary $2$-homomorphisms between those as
2-morphisms.
\end{proposition}

\noindent {\bf Proof. }  One need only check that if the
homomorphisms $F \maps L \to L'$ and $G \maps L' \to L''$ have
$F_2 = 1$, $G_2 = 1$, then their composite has $(FG)_2 = 1$.
\qed

The following proposition shows that strict Lie 2-algebras as
defined here agree with those as defined previously \cite{B}:

\begin{proposition} \et The 2-category {\rm SLie2Alg}
is isomorphic to the 2-category $\LieAlg\Cat$ consisting of
categories, functors and natural transformations in $\LieAlg$.
\end{proposition}

\noindent {\bf Proof. }  This is just a matter of unravelling the
definitions.  \qed

Just as Lie groups have Lie algebras, `strict Lie 2-groups' have
strict Lie 2-algebras.  Before we can state this result precisely,
we must recall the concept of a strict Lie $2$-group, which was
treated in greater detail in HDA5:

\begin{defn} \et \label{SLie2Grp} We define {\bf SLie2Grp} to be
the strict 2-category $\LieGrp\Cat$ consisting of categories,
functors and natural transformations in $\LieGrp$.  We call the
objects in this 2-category {\bf strict Lie 2-groups};
we call the morphisms between
these {\bf strict homomorphisms}, and we call the 2-morphisms
between those {\bf 2-homomorphisms}.
\end{defn}

\begin{proposition} \et  There exists a unique 2-functor
$$d \maps {\rm SLie2Grp} \to {\rm SLie2Alg}$$
such that:
\begin{enumerate}

\item $d$ maps any strict Lie $2$-group $C$ to the strict Lie $2$-algebra
$dC = \mathfrak{c}$ for which $\mathfrak{c}_0$ is the Lie algebra
of the Lie group of objects $C_0$, $\mathfrak{c}_1$ is the Lie
algebra of the Lie group of morphisms $C_1$, and the maps $s,t
\maps \mathfrak{c}_1 \rightarrow \mathfrak{c}_0$, $i \maps
\mathfrak{c}_0 \rightarrow \mathfrak{c}_1$ and $\circ \maps
\mathfrak{c}_1 \times_{\mathfrak{c}_0} \mathfrak{c}_1 \rightarrow
\mathfrak{c}_1$ are the differentials of those for $C$.

\item $d$ maps any strict Lie 2-group homomorphism
$F \maps C \rightarrow C'$ to the strict Lie 2-algebra
homomorphism $dF\maps \mathfrak{c} \rightarrow \mathfrak{c}'$ for
which $(dF)_0$ is the differential of $F_0$ and $(dF)_1$ is the
differential of $F_1$.

\item $d$ maps any strict Lie 2-group 2-homomorphism $\alpha \maps F \To G$
where $F,G \maps C \to C'$ to the strict Lie 2-algebra
2-homomorphism $d\alpha \maps dF \To dG$ for which the map
$d\alpha \maps \mathfrak{c}_0 \to \mathfrak{c}_1$ is the
differential of $\alpha \maps C_0 \to C_1$.
\end{enumerate}
\end{proposition}

\noindent{\bf Proof. } The proof of this long-winded proposition is a
quick exercise in internal category theory: the well-known functor
from $\LieGrp$ to $\LieAlg$ preserves pullbacks, so it maps
categories, functors and natural transformations in $\LieGrp$ to
those in $\LieAlg$, defining a 2-functor $d \maps {\rm SLie2Grp}
\to {\rm SLie2Alg}$, which is given explicitly as above.  \qed

We would like to generalize this theorem and define the Lie
2-algebra not just of a strict Lie 2-group, but of a general Lie
2-group as defined in HDA5.  However, this may require a weaker
concept of Lie 2-algebra than that studied here.


A nice way to obtain strict Lie 2-algebras is from `differential
crossed modules'.  This construction resembles one in HDA5, where 
we obtained strict Lie $2$-groups from `Lie crossed modules'.  
We recall that construction here before stating its Lie algebra analogue.

Starting with a strict Lie $2$-group $C$, with Lie groups
$C_{0}$ of objects and $C_{1}$ of morphisms, we construct a pair
of Lie groups
$$G=C_{0}, \; \; \; \; H = \ker(s) \subseteq C_{1}$$
where $s \colon C_1 \to C_0$ is the source map.
We then restrict the target map to a homomorphism
$$t \colon H \rightarrow G.$$  
In addition to the usual action of $G$ on
itself by conjugation, we have an action of $G$ on $H$,
$$\alpha \colon G \rightarrow \mathrm{Aut}(H),$$
defined by
$$\alpha(g)(h) = i(g)hi(g)^{-1}.$$
where $i \colon C_{0} \rightarrow C_{1}$ is the identity-assigning
map.  The target map is equivariant with respect to
this action, meaning:
$$t(\alpha(g))(h) = gt(h)g^{-1}.$$
We also have the `Peiffer identity':
$$\alpha(t(h))(h') = hh'h^{-1}$$
for all $h, h' \in H$.  So, we obtain the Lie group version of 
a crossed module:

\begin{defn} \et
A {\bf Lie crossed module} is a quadruple $(G,H,t,\alpha)$
consisting of Lie groups $G$ and $H$, a homomorphism $t: H
\rightarrow G$, and an action $\alpha$ of $G$ on $H$ (that is, a
homomorphism $\alpha \colon G \to {\rm Aut}(H)$) satisfying
$$t(\alpha(g)(h)) = g \, t(h)\, g^{-1} $$ and
$$   \alpha(t(h))(h') = hh'h^{-1}  $$
for all $g \in G$ and $h,h' \in H$.
\end{defn}

\noindent
In Proposition 32 of HDA5 we sketched how one can reconstruct
a strict Lie 2-group from its Lie crossed module.

For a Lie algebra analogue of this result, we should 
replace the Lie group $\mathrm{Aut}(H)$ by the Lie algebra
$\mathfrak{der}(\mathfrak{h})$ consisting of all `derivations' of
$\mathfrak{h}$, that is, all linear maps $f \colon \mathfrak{h}
\to \mathfrak{h}$ such that
$$  f([y,y']) = [f(y), y'] + [y, f(y')].  $$

\begin{defn} \et A {\bf differential crossed module} is a
quadruple $(\mathfrak{g}, \mathfrak{h}, t, \alpha)$ consisting of
Lie algebras $\mathfrak{g}$ and $\mathfrak{h}$, a homomorphism $t \colon
\mathfrak{h} \rightarrow \mathfrak{g}$, and an action $\alpha$ of
$\mathfrak{g}$ as derivations of $\mathfrak{h}$ (that is, a
homomorphism $\alpha: \mathfrak{g} \rightarrow
\mathfrak{der}(\mathfrak{h})$) satisfying 
$$  \alpha(x)(y) = [x,t(y)]  $$
and
$$  \alpha(t(y))(y') = [y,y']  $$
for all $x \in \mathfrak{g}$ and $y, y' \in \mathfrak{h}$.
\end{defn}

Differential crossed modules first appeared in the work of
Gerstenhaber \cite{G} where he classified them using the 3rd Lie algebra
cohomology group of $\mathfrak{g}$.   
We shall see a similar classification of Lie 
2-algebras in Corollary \ref{class}.  Indeed, differential
crossed modules are essentially the same as strict Lie 2-algebras:

\begin{proposition} \et Given a strict Lie 2-algebra $\mathfrak{c}$, 
there is a differential crossed module 
$(\mathfrak{g}, \mathfrak{h}, t, \alpha)$ where
$\mathfrak{g} = \mathfrak{c}_0$, $\mathfrak{h} = \ker(s)$, 
$t \colon \mathfrak{h} \to \mathfrak{g}$ is the restriction of 
the target map from $\mathfrak{c}_1$ to $\mathfrak{h}$, and 
$$  \alpha(x)(y) = [1_x, y]. $$ 
Conversely, we can reconstruct any strict Lie 2-algebra up to
isomorphism from its differential crossed module.
\end{proposition}

\noindent{\bf Proof. } The proof is analogous to the standard
one relating 2-groups and crossed modules \cite{BLau,FB}.
\qed

\noindent
The diligent reader can improve this proposition by
defining morphisms and $2$-morphisms for differential crossed 
modules, and showing this gives a $2$-category equivalent
to the $2$-category SLie2Alg.

Numerous examples of Lie crossed modules are described in Section
8.4 of HDA5.  Differentiating them gives examples of differential
crossed modules, and hence strict Lie $2$-algebras.



\section{Skeletal Lie $2$-algebras} \label{skeletallie2algs}

A semistrict Lie 2-algebra is {\it strict}
when we assume the map $l_3$ vanishes in the corresponding
$L_\infty$-algebra, since this forces
the Jacobiator to be the identity.    We now investigate
the consequences of assuming the differential $d$ vanishes
in the corresponding $L_\infty$-algebra.  Thanks to the
formula
$$   d\vec{f} =  t(f) - s(f)   , $$
this implies that the source of any morphism in the Lie 2-algebra
equals its target.  In other words, the Lie 2-algebra is
{\it skeletal:}

\begin{defn} \et
A category is {\bf skeletal} if isomorphic objects are equal.
\end{defn}

Skeletal categories are useful in category theory
because every category is equivalent to a skeletal one
formed by choosing one representative of each isomorphism
class of objects \cite{Mac}.  The same sort of thing is
true in the context of 2-vector spaces:

\begin{lemma}\label{skeletal2vs}\et
Any $2$-vector space is equivalent, as an
object of $2\Vect$, to a skeletal one.
\end{lemma}

\noindent {\bf Proof:} Using the result of Theorem
\ref{equivof2vs} we may treat our 2-vector spaces as 2-term chain
complexes.  In particular, a 2-vector space is skeletal if the
corresponding 2-term chain complex has vanishing differential, and
two 2-vector spaces are equivalent if the corresponding 2-term
chain complexes are chain homotopy equivalent. So, it suffices to
show that any 2-term chain complex is chain homotopy equivalent to
one with vanishing differential. This is well-known, but the basic
idea is as follows. Given a 2-term chain complex
$$\xymatrix{
   C_{1}
   \ar[rr]^{d}
    && C_{0}}
$$
we express the vector spaces $C_{0}$ and $C_{1}$ as
$C_{0} = \im(d) \oplus C_{0}'$ and $C_{1} = \ker(d) \oplus X$
where $X$ is a vector space complement to $\ker(d)$ in $C_{1}$.
This allows us to define a $2$-term chain complex $C'$ with vanishing
differential:
$$\xymatrix{
   C_{1}'= \ker(d)
   \ar[rr]^{0}
   && C_{0}'}
.$$
The inclusion of $C'$ in $C$ can easily be extended to a
chain homotopy equivalence. {\hbox{\hskip 30em} \qed}

Using this fact we obtain a result that will ultimately
allow us to classify Lie 2-algebras:

\begin{proposition} \et \label{skeletal}
Every Lie 2-algebra is equivalent, as an object of
{\rm Lie2Alg}, to a skeletal one.
\end{proposition}

\noindent {\bf Proof:} Given a Lie 2-algebra $L$ we may
use Lemma \ref{skeletal2vs} to find an equivalence between
the underlying 2-vector space of $L$ and a skeletal
2-vector space $L'$.  We may then use this to transport the
Lie 2-algebra structure from $L$ to $L'$, and obtain
an equivalence of Lie 2-algebras between $L$ and $L'$.
\qed

It is interesting to observe that
a skeletal Lie 2-algebra that is also strict
amounts to nothing but a Lie algebra $L_{0}$ together
with a representation of $L_{0}$ on a vector space $L_{1}$.
This is the infinitesimal analogue
of how a strict skeletal 2-group $G$
consists of a group $G_0$ together with an action of $G_0$
as automorphisms of an abelian group $G_1$.  Thus,
the representation theory of groups and Lie algebras is
automatically subsumed in the theory of 2-groups and Lie
2-algebras!

To generalize this observation to other skeletal Lie
2-algebras, we recall some basic definitions concerning
Lie algebra cohomology:


\begin{defn} \et
Let $\mathfrak{g}$ be a Lie algebra and $\rho$ a representation
of $\mathfrak{g}$ on the vector space $V$.  Then a
\textbf{V-{\bf valued} n-{\bf cochain}} $\mathbf{\omega}$ on
$\mathfrak{g}$ is a totally antisymmetric map
$$\omega \maps \mathfrak{g}^{\otimes n} \to V.$$
The vector space of all $n$-cochains is denoted
by $C^{n}(\mathfrak{g},V)$.
The {\bf coboundary operator} $\delta\maps C^{n}(\mathfrak{g},V)
\rightarrow C^{n+1}(\mathfrak{g},V)$ is defined by:
\begin{eqnarray*}
 (\delta \omega) (v_{1}, v_{2}, \dots, v_{n+1}) &:=&
    \sum_{i=1} ^{n+1} (-1)^{i+1} \rho (v_{i})
    \omega_{n} (v_{1}, \dots, \hat{v_{i}}, \dots, v_{n+1}) \cr
 &+& \sum _{1 \leq j < k \leq n+1} (-1)^{j+k} \omega _{n}
    ([v_{j}, v_{k}], v_{1}, \dots, \hat{v_{j}}, \dots,
    \hat{v_{k}}, \dots, v_{n+1}),
\end{eqnarray*}
\end{defn}

\begin{proposition} \et The Lie algebra coboundary operator $\delta$
satisfies $\delta ^{2} = 0$.
\end{proposition}

\begin{defn} \et
A $V$-valued $n$-cochain $\omega$ on $\mathfrak{g}$
is an {\bf $n$-cocycle} when $\delta
\omega = 0$ and an {\bf $n$-coboundary} if there exists an
$(n-1)$-cochain $\theta$ such that
$\omega = \delta \theta.$  We denote the groups of
$n$-cocycles and $n$-coboundaries
by $Z^{n}(\mathfrak{g},V)$ and $B^{n}(\mathfrak{g},V)$
respectively. The $n$th {\bf Lie algebra cohomology group}
$H^{n}(\mathfrak{g},V)$ is defined by
$$H^{n}(\mathfrak{g},V) = Z^{n}(\mathfrak{g},V)/B^{n}(\mathfrak{g},V).$$
\end{defn}

The following result illuminates the relationship between Lie
algebra cohomology and $L_{\infty}$-algebras.

\begin{theorem} \et \label{trivd}
There is a one-to-one correspondence between $L_{\infty}$-algebras
consisting of only two nonzero terms $V_{0}$ and $V_{n}$, with
$d=0,$ and quadruples $(\mathfrak{g}, V, \rho, l_{n+2})$ where
$\mathfrak{g}$ is a Lie algebra, $V$ is a vector space, $\rho$ is
a representation of $\mathfrak{g}$ on $V$, and $l_{n+2}$ is a
$(n+2)$-cocycle on $\mathfrak{g}$ with values in $V$.
\end{theorem}

\noindent {\bf Proof. }

\noindent $(\Rightarrow)$ Given such an $L_{\infty}$-algebra $V$
we set $\mathfrak{g}= V_{0}$.
$V_0$ comes equipped with a bracket as part of the
$L_{\infty}$-structure, and since $d$ is trivial, this bracket
satisfies the Jacobi identity on the nose, making
$\mathfrak{g}$ into a Lie algebra. We define $V = V_{n},$
and note that the bracket also gives a map
$\rho\maps \mathfrak{g} \otimes V \rightarrow V$, defined by
$\rho(x)f = [x,f]$ for $x \in \mathfrak{g}, f \in V$. We have
\begin{eqnarray*}
  \rho ([x,y])f &=& [[x,y],f] \cr
                &=& [[x,f],y] + [x,[y,f]] \; \; \; \;
                    \textrm{by $(2)$ of Definition \ref{L-alg}} \cr
                &=& [\rho(x)f, y] + [x, \rho(y) f]
\end{eqnarray*}

\noindent
for all $x,y \in \mathfrak{g}$ and $f \in V$, so that $\rho$ is
a representation.  Finally, the $L_{\infty}$ structure
gives a map
$l_{n+2}\maps \mathfrak{g}^{\otimes(n+2)} \rightarrow V$
which is in fact a $(n+2)$-cocycle.
To see this, note that
$$0 = \sum_{i+j = n+4} \sum_{\sigma}
l_{j}(l_{i}(x_{\sigma(1)}, \ldots, x_{\sigma(i)}),
x_{\sigma(i+1)}, \ldots, x_{\sigma(n+2)})$$
where we sum over $(i, (n+3)-i)$-unshuffles $\sigma \in S_{n+3}$.
However, the only choices
for $i$ and $j$ that lead to nonzero $l_{i}$ and $l_{j}$ are
$i=n+2, j=2$ and $i=2, j=n+2.$  In addition, notice that in this
situation, $\chi(\sigma)$ will consist solely of the sign of the
permutation because all of our $x_{i}$'s have degree zero.  Thus,
the above becomes, with $\sigma$ a $(n+2, 1)$-unshuffle and $\tau$
a $(2, n+1)$-unshuffle:

\begin{eqnarray*}
  0 &=& \sum_{\sigma} \chi(\sigma) (-1)^{n+2}
     [l_{n+2}(x_{\sigma(1)}, \dots, x_{\sigma(n+2)}), x_{\sigma(n+3)}] \\
     &&+ \sum_{\tau} \chi(\tau) l_{n+2}([x_{\tau(1)}, x_{\tau(2)}],
     x_{\tau(3)}, \dots, x_{\tau(n+3)}) \cr  \\
    &=& \sum_{i=1}^{n+3} (-1)^{n+3-i}(-1)^{n+2}
     [l_{n+2}(x_{1}, \dots, x_{i-1}, x_{i+1},\dots, x_{n+3}), x_{i}] \cr  \\
    & & + \sum_{1 \leq i < j \leq n+3} (-1)^{i+j+1}
     l_{n+2}([x_{i}, x_{j}], x_{1}, \dots,\hat{x_{i}},
     \dots, \hat{x_{j}}, \dots, x_{n+3})
     \qquad \qquad (\dag) \cr \\
    &=& \sum_{i=1}^{n+3} (-1)^{i+1} [l_{n+2}(x_{1}, \dots, x_{i-1},
     x_{i+1}, \dots, x_{n+3}), x_{i}] \cr  \\
    & & + \sum_{1 \leq i < j \leq n+3} (-1)^{i+j+1}
     l_{n+2}([x_{i}, x_{j}], x_{1}, \dots, \hat{x_{i}},
    \dots, \hat{x_{j}}, \dots, x_{n+3}) \cr \\
    &=& - \sum_{i=1}^{n+3} (-1)^{i+1} [x_{i},
     l_{n+2}(x_{1}, \dots, x_{i-1}, x_{i+1}, \dots, x_{n+3})]  \cr  \\
    & &  - \sum_{1 \leq i < j \leq n+3} (-1)^{i+j}
     l_{n+2}([x_{i}, x_{j}], x_{1}, \dots, \hat{x_{i}}, \dots,
      \hat{x_{j}}, \dots, x_{n+3})  \cr \\
    &=& -\delta l_{n+2}(x_{1}, x_{2}, \dots, x_{n+3}).
\end{eqnarray*}

\noindent The first sum in $(\dag)$ follows because we have
$(n+3)$ $(n+2,1)$-unshuffles and the sign of any such unshuffle is
$(-1)^{n+3-i}$. The second sum follows similarly because we have
$(n+3)$ $(2, n+1)$-unshuffles and the sign of a
$(2,n+1)$-unshuffle is $(-1)^{i+j+1}.$  Therefore, $l_{n+2}$ is a
$(n+2)$-cocycle.

\vskip 1em \noindent $(\Leftarrow)$ Conversely, given a Lie
algebra $\mathfrak{g}$, a representation $\rho$ of $\mathfrak{g}$
on a vector space $V$, and an $(n+2)$-cocycle $l_{n+2}$ on
$\mathfrak{g}$ with values in $V$, we define our
$L_{\infty}$-algebra $V$ by setting
$V_{0} = \mathfrak{g}$, $V_{n} = V$, $V_i = \{0\}$ for $i \ne 0,n$,
and $d=0.$  It remains to define the system of linear maps
$l_{k}$, which we do as follows: Since $\mathfrak{g}$ is a Lie
algebra, we have a bracket defined on $V_{0}$. We extend this
bracket to define the map $l_2$, denoted by $[\cdot, \cdot] \maps
V_{i} \otimes V_{j} \rightarrow V_{i+j}$ where $i,j=0,n,$ as
follows:
$$[x,f] = \rho (x) f,$$
$$[f,y] = - \rho(y)f,$$
$$[f,g] =0$$
for $x,y \in V_0$ and $f,g \in V_n$.
With this
definition, the map $[\cdot, \cdot]$ satisfies condition $(1)$ of
Definition \ref{L-alg} . We define $l_{k}=0$ for $3 \leq k \leq
n+1$ and $k> n+2$, and take $l_{n+2}$ to be the given $(n+2)$ cocycle,
which satisfies conditions $(1)$ and $(2)$ of Definition \ref{L-alg}
by the cocycle condition. \qed

We can classify skeletal Lie $2$-algebras
using the above construction with $n=1$:

\begin{corollary} \et \label{class} There is a one-to-one correspondence
between isomorphism classes
of skeletal Lie $2$-algebras and isomorphism classes
of quadruples consisting of a Lie algebra $\mathfrak{g}$, a vector
space $V$, a representation $\rho$ of $\mathfrak{g}$ on $V$, and
a 3-cocycle on $\mathfrak{g}$ with values in $V$.
\end{corollary}

\noindent {\bf Proof. }  This is immediate from Theorem
\ref{1-1} and Theorem \ref{trivd}.  \qed

Since every Lie 2-algebra is equivalent as an object of
{\rm Lie2Alg} to a skeletal one,
this in turn lets us classify {\it all} Lie 2-algebras,
though only up to equivalence:

\begin{theorem} \et \label{class2} There is a one-to-one correspondence
between equivalence classes
of Lie $2$-algebras (where equivalence is as objects
of the 2-category {\rm Lie2Alg}) and isomorphism classes
of quadruples consisting of a Lie algebra $\mathfrak{g}$, a vector
space $V$, a representation $\rho$ of $\mathfrak{g}$ on $V$, and
an element of $H^3(\mathfrak{g},V)$.
\end{theorem}

\noindent {\bf Proof. }  This follows from Theorem \ref{skeletal}
and Corollary \ref{class}; we leave it to the reader to verify
that equivalent skeletal Lie 2-algebras give
cohomologous 3-cocycles and conversely.  \qed

We conclude with perhaps the most interesting examples of
finite-dimensional Lie 2-algebras
coming from Theorem \ref{class}.  These make use of
the following identities involving the Killing form
$\langle x,y\rangle := \tr(\ad(x)\ad(y)) $
of a finite-dimensional Lie algebra:
$$\langle x, y \rangle = \langle y, x \rangle,$$
and
$$\langle [x,y], z \rangle = \langle x, [y,z] \rangle.$$

\begin{example} \label{ghbar} \et There is a skeletal Lie $2$-algebra
built using Theorem \ref{class} by taking
$V_{0} = \mathfrak{g}$ to be a finite-dimensional Lie algebra
over the field $k$, $V_{1}$ to be $k$, $\rho$ the
trivial representation, and
$l_3(x,y,z) = \langle x,[y,z]\rangle $.
We see that $l_3$ is a $3$-cocycle using
the above identities as follows:
\begin{eqnarray*}
  (\delta l_3)(w,x,y,z) &=& \rho(w) l_3(x,y,z) - \rho(x)l_3(w,y,z) +
                          \rho(y)l_3(w,x,z) - \rho(z) l_3(w,x,y) \cr
           & & -l_3([w,x],y,z) + l_3([w,y],x,z) -l_3([w,z],x,y) \cr
           & & -l_3([x,y],w,z) + l_3([x,z],w,y) -l_3([y,z],w,x) \cr
                      &=& -\langle [w,x],[y,z] \rangle
                          +\langle [w,y],[x,z] \rangle
                          -\langle [w,z],[x,y] \rangle \cr
                      & & -\langle [x,y],[w,z] \rangle
                          +\langle [x,z],[w,y] \rangle
                          -\langle [y,z],[w,x] \rangle
\end{eqnarray*}

\noindent This second step above follows because we have a trivial
representation. Continuing on, we have
\begin{eqnarray*}
   (\delta l_3)(w,x,y,z)  &=& -2\langle [w,x],[y,z] \rangle
                            +2\langle [w,y],[x,z] \rangle
                            -2\langle [w,z],[x,y] \rangle \cr
                        &=& -2\langle w,[x,[y,z]] \rangle
                            +2\langle w,[y,[x,z]] \rangle
                            -2\langle w,[z,[x,y]] \rangle \cr
                        &=& -2\langle w,[x,[y,z]] + [y,[z,x]]
                                +[z,[x,y]]\rangle  \cr
                        &=& -2\langle w, 0 \rangle   \cr
                        &=& 0. \; \; \;
\end{eqnarray*}
More generally, we obtain a Lie 2-algebra this way taking
$l_3(x,y,z) = \hbar \langle x,[y,z]\rangle $
where $\hbar$ is any element of $k$.
We call this Lie 2-algebra $\mathfrak{g}_\hbar$.
\end{example}

It is well known that the Killing form of $\mathfrak{g}$ is
nondegenerate if and only if $\mathfrak{g}$ is semisimple.
In this case the 3-cocycle described above represents a
nontrivial cohomology class when $\hbar \ne 0$,
so by Theorem \ref{class2} the Lie 2-algebra $\mathfrak{g}_\hbar$
is not equivalent to a skeletal one with vanishing Jacobiator.
In other words, we obtain a Lie 2-algebra that is not equivalent
to a skeletal strict one.

Suppose the field $k$ has characteristic zero, the
Lie algebra $\mathfrak{g}$ is finite dimensional and
semisimple, and $V$ is finite dimensional.  Then a
version of Whitehead's Lemma \cite{AIP} says that
$H^3(\mathfrak{g},V) = \{0\}$ whenever the representation
of $\mathfrak{g}$ on $V$ is nontrivial and irreducible.
This places some limitations on finding interesting
examples of nonstrict Lie 2-algebras other than those
of the form $\mathfrak{g}_\hbar$.

In HDA5 we show how the Lie 2-algebras $\g_\hbar$ give rise
to 2-groups when $\hbar$ is an integer.  The construction
involves Chern--Simons theory.  Since Chern--Simons theory
is also connected to the theory of quantum groups and affine
Lie algebras, it is natural to hope for a more direct link
between these structures and the Lie 2-algebras $\g_\hbar$.
After all, they are all `deformations' of more familiar
algebraic structures which take advantage of the
3-cocycle $\langle x, [y,z] \rangle$ and the closely
related 2-cocycle on $C^\infty(S^1,\mathfrak{g})$.

The smallest nontrivial example of the Lie 2-algebras
$\g_\hbar$ comes from $\g = \mathfrak{su}(2)$.
Since $\mathfrak{su}(2)$ is isomorphic to $\R^3$ with
its usual vector cross product, and its Killing form is
proportional to the dot product, this Lie 2-algebra relies
solely on familiar properties of the dot product and
cross product: $$x \times y = - y \times x,$$
$$x \cdot y = y \cdot x,$$
$$x \cdot (y \times z) = (x \times y) \cdot z,$$
$$x \times (y \times z) + y \times (z \times x) + z \times (x \times y) =0.$$
It will be interesting to see if this Lie 2-algebra, where
the Jacobiator comes from the triple product, has any
applications to physics.
Just for fun, we work out the details again in this case:

\begin{example} \et There is a skeletal Lie $2$-algebra
built using Theorem \ref{class} by taking
$V_{0} = \R^3$ equipped with the cross product, $V_{1} = \R$,
$\rho$ the trivial representation, and $l_3(x,y,z) = x \cdot (y \times z)$.
We see that $l_3$ is a 3-cocycle as follows:

\begin{eqnarray*}
   (\delta l_3)(w,x,y,z) &=&
      -l_3([w,x],y,z) + l_3([w,y],x,z) -l_3([w,z],x,y) \cr
                & & -l_3([x,y],w,z) + l_3([x,z],w,y) -l_3([y,z],w,x) \cr
                       &=& -(w \times x) \cdot (y \times z)
                           +(w \times y) \cdot (x \times z)
                           -(w \times z) \cdot (x \times y) \cr
                       & & -(x \times y) \cdot (w \times z)
                           +(x \times z) \cdot (w \times y)
                           -(y \times z) \cdot (w \times x) \cr
                       &=& -2(w \times x) \cdot (y \times z)
                           +2(w \times y) \cdot (x \times z)
                           -2(w \times z) \cdot (x \times y) \cr
                       &=& -2w \cdot (x \times (y \times z))
                           +2w \cdot (y \times (x \times z))
                           -2w \cdot (z \times (x \times y)) \cr
                       &=& -2w \cdot (x \times (y \times z)
                     +y \times (z \times x) + z \times (x \times y)) \cr
                      &=& 0.
\end{eqnarray*}

\end{example}


\section{Conclusions} \label{conclusions}

In HDA5 and the present paper we have seen evidence that the theory
of Lie groups and Lie algebras can be categorified to give interesting
theories of Lie 2-groups and Lie 2-algebras.  We expect this pattern
to continue as shown in the following tables.

\vskip 0.5em
\begin{center}
{\small
\begin{tabular}{|c|c|c|c|}  \hline
      & $n = 0$   & $n = 1$    & $n = 2$                  \\     \hline
$k = 0$  & manifolds  & Lie groupoids  & Lie 2-groupoids  \\     \hline
$k = 1$  & Lie groups & Lie 2-groups   & Lie 3-groups     \\     \hline
$k = 2$  & abelian    & braided        & braided          \\
         & Lie groups & Lie 2-groups   & Lie 3-groups     \\     \hline
$k = 3$  &`'          & symmetric      & sylleptic        \\
         &            & Lie 2-groups   & Lie 3-groups     \\     \hline
$k = 4$  &`'          &`'              & symmetric        \\
         &            &                & Lie 3-groups     \\     \hline
$k = 5$  &`'          &`'              & `'               \\
         &            &                &                  \\     \hline
\end{tabular}} \vskip 1em
1. $k$-tuply groupal Lie $n$-groupoids:
expected results
\end{center}

\vskip 0.5em
\begin{center}
{\small
\begin{tabular}{|c|c|c|c|}  \hline
      & $n = 0$   & $n = 1$    & $n = 2$                  \\     \hline
$k = 0$  & vector bundles & Lie algebroids & Lie 2-algebroids \\     \hline
$k = 1$  & Lie algebras   & Lie 2-algebras & Lie 3-algebras   \\     \hline
$k = 2$  & abelian        & braided        & braided          \\
         & Lie algebras   & Lie 2-algebras & Lie 3-algebras   \\     \hline
$k = 3$  &`'              & symmetric      & sylleptic        \\
         &                & Lie 2-algebras & Lie 3-algebras   \\     \hline
$k = 4$  &`'              & `'             & symmetric        \\
         &                &                & Lie 3-algebras     \\     \hline
$k = 5$  &`'              &`'              & `'               \\
         &                &                &                  \\     \hline
\end{tabular}} \vskip 1em
2. $k$-tuply stabilized Lie $n$-algebroids:
expected results
\end{center}
\vskip 0.5em

Table 1 gives names for $k$-tuply groupal $n$-groupoids
\cite{BD} for which the set of $j$-morphisms is a smooth
manifold for each $j$, and for which the operations are all
smooth.  Manifolds, Lie groups and abelian Lie groups are well-understood;
Lie groupoids have also been intensively investigated \cite{Mackenzie},
but the study of Lie 2-groups has just barely begun, and the
other entries in the chart are still {\it terra incognita}:
they seem not to have even been defined yet, although this should
be easy for the entries in the second column.

Table 2 gives names for the `infinitesimal versions'
of the entries in the first chart.  The classic example is that
of a Lie algebra, which can be formed by taking the tangent space
of a Lie group at the identity element.  Similarly, we have
seen that the tangent 2-vector space at the identity object of
a strict Lie 2-group becomes a strict Lie 2-algebra; we also
expect a version of this result to hold for the more general Lie
2-groups defined in HDA5, though our `semistrict Lie 2-algebras'
may not be sufficiently general for this task.

The $k = 0$ row of Table 2 is a bit different from the rest.
For example, a manifold does not have a distinguished identity
element at which to take the tangent space.   To deal with
this we could work instead with pointed manifolds, but another
option is to take the tangent space at {\it every} point of a manifold
and form the tangent bundle, which is a vector bundle.  Similarly,
a Lie groupoid does not have a distinguished `identity object',
so the concept of `Lie algebroid' \cite{Mackenzie} must be
defined a bit subtly.   The same will be true of Lie $n$-groupoids
and their Lie $n$-algebroids.  For this reason it may be useful
to treat the $k = 0$ row separately and use the term `$k$-tuply
stabilized Lie $n$-algebra' for what we are calling a $(k+1)$-tuply
stabilized Lie $n$-algebroid.

The general notion of `$k$-tuply stabilized Lie $n$-algebra'
has not yet been defined, but at least we understand
the `semistrict' ones: as explained in Section \ref{Linftyalgs},
these are just various sorts of $L_\infty$-algebra with their
underlying $n$-term chain complexes reinterpreted as strict
$(n-1)$-categories in $\Vect$.  More precisely, we define a
{\bf semistrict $k$-tuply stabilized Lie $n$-algebra} to be the
result of taking an $L_\infty$-algebra $V$ with $V_i = 0$ when
$i < k$ or $i \ge n + k$ and transferring all the structure on
its underlying $n$-term chain complex to the corresponding
strict $(n-1)$-category in $\Vect$.

In this language, Theorem \ref{trivd} gives a way of constructing
a semistrict Lie $n$-algebra with only nontrivial objects
and $n$-morphisms from an $(n+2)$-cocycle on the Lie algebra
of objects.  This can be seen as an infinitesimal version of the usual
`Postnikov tower' construction of a connected homotopy $(n+1)$-type
with only $\pi_1$ and $\pi_{n+1}$ nonzero from an $(n+2)$-cocycle
on the group $\pi_1$.  The analogy comes into crisper focus
if we think of a connected homotopy $(n+1)$-type as an `$n$-group'.
Then the Postnikov construction gives an $n$-group with only
nontrivial objects and $n$-morphisms from an $(n+2)$-cocycle on
the group of objects; now the numbering scheme perfectly matches that
for Lie $n$-algebras.   For $n = 2$ we described how this works
more explicitly in HDA5.   One of the goals of the present paper
was to show that just as group cohomology arises naturally in the
classification of $n$-groups, Lie algebra cohomology arises
in the classification of Lie $n$-algebras.


\subsection*{Acknowledgements}

We thank Ronnie Brown, James Dolan, Vyjayanthi Chari, Andr\'ee
Ehresmann, Aaron Lauda, Thomas Larsson and James Stasheff for
helpful discussions and correspondence.  We also thank J.\ Scott
Carter and Masahico Saito for letting us use pictures of the
Zamolodchikov tetrahedron equation taken from their book \cite{CS},
and even adapting these pictures specially for us.  Finally, we thank Behrang Noohi \cite{Noohi} and Dmitry Roytenberg \cite{Roytenberg} for helping us catch and fix some significant errors in this paper.


\end{document}